\documentclass[10pt]{article}
\usepackage{CJK}
\usepackage{mathtools}
\usepackage{cases}
\usepackage{graphicx,amssymb,mathrsfs,amsmath,color}
\usepackage{pifont}
\newtheorem{definition}{Definition}

\newtheorem{lemma}{Lemma}
\newtheorem{theorem}{Theorem}
\newtheorem{corollary}{Corollary}[section]

\newcommand{\C}[1]{{\cal {#1}}}
\newcommand{\tr}{^{\sf T}}
\newcommand{\m}[1]{{\bf{#1}}}

\begin{document}
\title{Iteration complexity  analysis of  a partial  LQP-based  alternating direction method of multipliers
\thanks{Authors are ranked by family name with equal contributions.
Part of this research was   supported by    National  Natural Science Foundation of China (12001430, 72071158),   Fundamental Research Funds for the Central Universities (G2020KY05203),  and   China Postdoctoral Science Foundation (2020M683545).}}
\author{
   Jianchao Bai\thanks{ {\tt jianchaobai@nwpu.edu.cn}, https://teacher.nwpu.edu.cn/jcbai,
        School of Mathematics and Statistics  and the MIIT Key Laboratory of Dynamics and Control of Complex Systems,  Northwestern Polytechnical
        University, Xi'an  710129,       China.}
\and
  Yuxue Ma\thanks{ \ding {41}  Corresponding author. {\tt  Mayuxue708991@163.com},
        School of Mathematics and Statistics,  Northwestern Polytechnical
        University, Xi'an  710129,       China.}
\and
Hao Sun \thanks{ \ding {41} Corresponding author. {\tt  hsun@nwpu.edu.cn},
        School of Mathematics and Statistics,
        Northwestern Polytechnical University, Xi'an  710129,      China.}
\and
Miao Zhang \thanks{ {\tt mzhan33@lsu.edu},
        Department of Mathematics, Louisiana State University, Baton Rouge, LA 70803-4918, USA}
}

\date{ }
\maketitle

\begin{abstract}
In this paper, we consider a prototypical  convex optimization problem with multi-block variables and separable structures.  By adding the Logarithmic Quadratic Proximal (LQP)   regularizer with suitable proximal parameter to each of the first grouped subproblems, we develop a partial  LQP-based   Alternating Direction Method of Multipliers  (ADMM-LQP).  The dual variable   is updated twice with relatively larger stepsizes than the classical region $(0,\frac{1+\sqrt{5}}{2})$. Using   a prediction-correction approach to analyze properties of the iterates generated by ADMM-LQP, we  establish its global convergence and   sublinear convergence rate of $O(1/T)$ in the new ergodic and nonergodic senses, where $T$ denotes the iteration index. We also extend the algorithm to a nonsmooth composite convex optimization and establish {similar   convergence results} as our ADMM-LQP.
\end{abstract}

\bigskip

\section{Introduction}

Consider the following   multi-block separable convex optimization problem
\begin{equation} \label{Sec1-prob}
\begin{array}{lll}
 \min & \sum\limits_{i=1}^{p}f_i(x_i)+g(y)\\
\textrm{s.t. }& \sum\limits_{i=1}^{p}A_ix_i+ By =b,\\
       & x_i \in\mathbb{R}_{+}^{m_i}{,}~y \in\C{Y}~,i=1,\cdots,p,
\end{array}
\end{equation}
where $f_i(x_i):\mathbb{R}_{+}^{m_i}\to\mathbb{R},\;g(y):\mathbb{R}^{d}\to\mathbb{R}$ are
closed  proper convex functions (but not necessarily  smooth/{ strongly convex); $A_i\in\mathbb{R}^{n\times m_i}, B\in\mathbb{R}^{n\times d}$ and $b\in \mathbb{R}^{n}$ are  given data;}   $\C{Y}\subset\mathbb{R}^{d}$ is a closed convex set and $p\ge1$ is an integer. Throughout our discussions, the solution set of the problem (\ref{Sec1-prob}) is assumed to be nonempty and all of matrices $A_i(i=1,\ldots,p)$ and $B$ have full column ranks.  {For convenience of analysis,} we denote $A=(A_{1},\ldots,A_{p}),\;\m{x}=({x_{1}\tr },\ldots,{x_{p}\tr })\tr,\;
\C{X}=\mathbb{R}_{+}^{m_1}\times\cdots\times\mathbb{R}_{+}^{m_p}$ and $\C{W}=\C{X}\times\C{Y}\times\mathbb{R}^{n}$.

{
Problems in the form of (\ref{Sec1-prob}) cover  a large number of practical applications, for example,   traffic assignment, network economics, game theoretic and cancer diagnostics
problems,   see  \cite{Fpang03,HeYuanLQP06,VDS2019,KiSt80,WuLi2019}. Here we   take three concrete motivating examples:
\begin{itemize}
  \item
  \textbf{Traffic network equilibrium problem}. Consider a directed traffic network   including 25
   nodes, 37 links, and 6 origin/destination (O/D) pairs.   In order to avoid  traffic
congestion, it is common to require the link flow to be bounded from above, that is, the path  flow $x$ should satisfy $x\in S:=\{x\in  \mathbb{R}^{55}| A\tr x\leq b, x\geq \m{0}\},$
where $A\in \mathbb{R}^{55\times 37}$   is the path-link incidence matrix and $b$ is the given link capacity vector.
Then, the traffic network equilibrium problem can be described to seek a path flow
pattern $x^*$ such that
\begin{equation} \label{Sec1-prob0}
x^*\geq \m{0},\quad \left\langle x-x^*, f(x^*)\right\rangle\geq0,\quad \forall x\in S,
\end{equation}
where $f(x^*)= At(A\tr x)- E\eta (E\tr x)$, $E\in\mathbb{R}^{55\times 6}$ is the path-O/D pair incidence
matrix, $t$ represents a given link travel cost vector  and $\eta$ denotes  the travel disutility determined by O/D pair. As discussed in \cite{HeYuanLQP06},
the traffic equilibrium problem (\ref{Sec1-prob0}) can be mathematically characterized as the  form (\ref{Sec31-001}),  i.e., an  equivalent form of Problem (\ref{Sec1-prob})  with  $\theta_2(y)\equiv 0, B=\m{I},p=1$ and $\C{Y}=\mathbb{R}_{+}^{55}$. Performance of  the two-block case of our proposed algorithm for solving this problem had been verified in  the experiments  \cite[Section 5]{WuLi2019}.

\item \textbf{Sparse signal processing problem}. As stated in \cite{xu2012}, many signal processing problem in medical imaging, computed tomography and so forth can be mathematically modeled as a large-scale linear equations $\C{A}\m{x}=b$.  In order to find  a sparse nonnegative solution to this equation with e.g. $\C{A}\in \mathbb{R}^{10000\times 5000}$ and $b\in \mathbb{R}^{10000}$,   an effective way is to split $\C{A}=[A_1,A_2,\cdots,A_{10}]$ and $\m{x}=(x_1,x_2, \cdots,x_{10})\tr$ where $x_i\in\mathbb{R}_{+}^{500}$. Hence, the problem can be converted to  the following small-scale minimization  problem:
    \[
\begin{array}{lll}
\min  &  \|x_1\|_1 + \|x_2\|_1+\cdots+\|x_{10}\|_1\\
  \textrm{s.t.}  &  A_1 x_1 +A_2 x_2 +\cdots +A_{10} x_{10} =b,~x_i\geq \m{0}, i=1,\cdots,10.
    \end{array}
\]
Clearly, this reformulated problem is a multi-block case of Problem (\ref{Sec1-prob}) with $p=9,f_i=\|x_i\|_1,g=\|x_{10}\|_1, B=A_{10}$ and $\C{Y}=\mathbb{R}_{+}^{500}$.
\item \textbf{Linear programming with box  constraints}. Linear programming   can date back to the time of Fourier  who published a method to solve problems with a system of linear inequalities, see e.g. \cite[page 257]{BoydVan10}. It is   a widely studied field in optimization and covers many practical problems, e.g. network flow problems and the model of trolley network \cite{LHLLZZ}. Linear programming problems with box-type constraints can be expressed as
\[
\begin{array}{lll}
 \min & c\tr z\\
\textrm{s.t. }& B z =b,l \leq z \leq u,
\end{array}
\]
where $c,   l, u\in \mathbb{R}^n$, $B\in \mathbb{R}^{m\times n}$  and $b\in\mathbb{R}^m$ are given. Its dual problem reads
\[
\begin{array}{lll}
 \min &       u\tr x_1- l\tr x_2  + b\tr y\\
\textrm{s.t. }&    x_1 -x_2   +B\tr y=-c, \\
&x_1\ge \m{0}, x_2 \ge \m{0}, y\in \mathbb{R}^m,
\end{array}
\]
which is clearly the form of Problem (\ref{Sec1-prob}) with three blocks.
\end{itemize}}

A prototype method for solving the equality constrained problem {in the form of}  (\ref{Sec1-prob}) is  the     Augmented Lagrangian
Method (ALM, \cite{Hestenes69}). For any $ \beta>0$, by constructing the augmented Lagrange function associated with  (\ref{Sec1-prob}):
\begin{equation} \label{Sec1-3}
\C{L}_\beta(\m{x},y,\lambda)= \C{L}(\m{x},y,\lambda)+\frac {\beta}{2}\left\|A \m{x}
+By-b\right\|^2,
 \end{equation}
where  \begin{equation} \label{Sec1-lagj}
\qquad~ \C{L}(\m{x},y,\lambda)= \sum\limits_{i=1}^{p}f_i(x_i)
+g(y)-\left\langle\lambda, A \m{x}
+By-b\right\rangle,
 \end{equation}
ALM   firstly updates  the primal variables  by solving the   jointed subproblem
\[
 \min _{(\m{x},y)\in\C{X}\times\C{Y}}  {\C{L}_\beta}(\m{x},y,\lambda^k)
\]
and then   updates the Lagrange multipliers by $\lambda^{k+1}=\lambda^k-\beta(A \m{x}^{k+1}
+By^{k+1}-b)$. However, ALM does not make full use of the separable structure of the
objective functions and hence, could not take advantage of {some} special
properties of each component function; furthermore, in many real applications involving big-data, solving this jointed subproblem
{would be} very expensive or even difficult at each iteration.

An effective approach to overcome the above disadvantages is the Alternating
Direction Method of Multipliers (ADMM)  which could be regarded as a
splitting version of ALM:
\begin{equation} \label{Sec1-2b}
\quad\left\{
  \begin{array}{lll}
    x_1^{k+1}&\leftarrow&\arg\min\limits_{x_1\in \mathbb{R}_{+}^{m_1}} \C{L}_{\beta}(x_1,x_2^k,\cdots,x_p^k,y^{k},\lambda^{k}) ,\\
    ~~\vdots&&\\
    x_p^{k+1}&\leftarrow&\arg\min\limits_{x_p\in \mathbb{R}_{+}^{m_p}} \C{L}_{\beta}(x_1^{k+1},\cdots,x_{p-1}^{k+1},x_p,y^{k},\lambda^{k}) ,\\
    y^{k+1}&\leftarrow&\arg\min\limits_{y\in \C{Y}} \C{L}_{\beta}(x_1^{k+1},\cdots,x_p^{k+1},y,\lambda^{k}) ,\\
    \lambda^{k+1}&\leftarrow&\lambda^k-\beta(A \m{x}^{k+1}
+By^{k+1}-b).
  \end{array}
\right.
 \end{equation}
{Compared  ADMM (\ref{Sec1-2b}) to ALM, each subproblem of ADMM  annihilates the coupled term $\beta\langle A\m{x}, By\rangle$  and makes its iteration independent on other variables,}   so ADMM can be quite effective for solving large-scale problems  with preferred structures. As commented by
Boyd et al. \cite{Boyd2010}, ``ADMM is at least comparable to very specialized algorithms (even in the
serial setting), and in most cases, the simple ADMM algorithm will be efficient enough to be
useful.''
Although the classical ADMM   \cite{DB1976,GlMa75} was  demonstrated convergent for the two-block separable convex optimization, its direct extension in the   Gauss-Seidel scheme (\ref{Sec1-2b}) is not necessarily convergent  \cite{chyy2016} without special assumptions on the coefficient matrix in the equality constraints. Besides, it was   pointed out by
  He et al. \cite{HeHouYuan2015} that the following Jacobi-type  extension of ADMM:
\begin{equation}  \label{Sec1-DE-ADMM2}
\left\{
  \begin{array}{lll}
    x_1^{k+1}&\leftarrow&\arg\min\limits_{x_1\in \mathbb{R}_{+}^{m_1}} \C{L}_{\beta}(x_1,x_2^k,\cdots,x_p^k,y^{k},\lambda^{k}) ,\\
    ~~\vdots&&\\
    x_p^{k+1}&\leftarrow&\arg\min\limits_{x_p\in \mathbb{R}_{+}^{m_p}} \C{L}_{\beta}(x_1^k,\cdots,x_{p-1}^k,x_p,y^{k},\lambda^{k}) ,\\
    y^{k+1}&\leftarrow&\arg\min\limits_{y\in \C{Y}} \C{L}_{\beta}(x_1^{k+1},\cdots,x_p^{k+1},y,\lambda^{k}) ,\\
    \lambda^{k+1}&\leftarrow&\lambda^k-\beta(A \m{x}^{k+1}
+By^{k+1}-b),
  \end{array}
\right.
\end{equation}
 is   not necessarily convergent either.

 Since  the development of ADMM in 1970s, a great majority of   ADMM-type methods  were proposed     { for solving}  the multi-block separable convex optimization problems. Particularly, many researchers focused on modifying the above scheme  (\ref{Sec1-2b}) or  (\ref{Sec1-DE-ADMM2}) by adding  proper
proximal terms into the subproblems of ADMM{, which aims  to promote  the algorithmic} convergence   and also to overcome the {obstacle that each involved subproblem is not easily solvable}. For instance, a proximal ADMM using the Gauss-Seidel approach  (\ref{Sec1-2b}) was developed in \cite{Hetalk}
to solve the $3$-block separable convex problem, in which proximal terms $\frac{\tau \beta}{2}\|A_i(\m{x}_i-\m{x}_i^k)\|^2$  were
added for the second and   third subproblems, and global  convergence was guaranteed when the proximal parameter  $\tau\in  [1,+\infty)$.
To deal with the convex problem (\ref{Sec1-prob}) with $N$-block variables ($N \ge 3$),  a  proximal ADMM using the Jacobian approach  (\ref{Sec1-DE-ADMM2}) was proposed in \cite{HeXuYuan2016b},
where  the   terms in the form {of} $\frac{s \beta}{2}\|A_i(\m{x}_i-\m{x}_i^k)\|^2$ were added to each generated subproblem and {its} global convergence
was established when the proximal parameter $s\in [N-1,+\infty)$.
 In {the recent work} \cite{He2015yuan,HeMaYu016}, the proximal point techniques were further studied and
  global convergence of {the   algorithms therein} was   established with the aid of a  prediction-correction approach for the iterative sequence.
More recently, by designing a  positive proximal  term  for each $x_i$-subproblem, Bai et al. \cite[the scheme (74)]{BaiLiXuZhang2018} developed  a generalized symmetric ADMM whose dual variable was updated twice with larger stepsizes than the   region $(0,\frac{1+\sqrt{5}}{2})$. Under the assumptions that all subdifferential  of each component objective function is piecewise linear multi-functions, Bai et al. \cite{BXCL191} even established the linear convergence {rate} of this generalized symmetric ADMM.
We refer   interested readers to \cite{CWHLv19,FPST14,ZhAD19,Li2012,XuY2017,WuLi2019} for other  proximal ADMMs by {exploiting   general  proximal terms
to solve  the two block and multiple block cases}.

It is a fact acknowledged that
the performance of ADMM-type methods   depends significantly on the difficulty of solving  each  subproblem. For the problem (\ref{Sec1-prob}), the involved   subproblems usually do
not have  closed-form solutions and   need  to be solved approximately by   {inexact methods or some   inner iterative algorithms. Fortunately, for the nonnegative linearly   constrained problem} the Logarithmic Quadratic Proximal (LQP, \cite{ATe99}) regularizer can convert the splitting subproblem to a system of nonlinear \textcolor{blue}{equations}, which is comparatively   easier
  than the original one. Another merit is that  LQP regularization can
ensure that the solution of involved subproblem stays strictly within the interior of positive
orthant. Mainly motivated by the advantages of  LQP regularization and  previous researches \cite{BaiLiXuZhang2018,liuGY19,WuLL2016,WuLi2019}, in this paper  we propose a partial  LQP-based   ADMM  (\textbf{ADMM-LQP}),  which is described formally in the following table,  with larger stepsizes of dual variables to solve the general  problem (\ref{Sec1-prob}).  In our proposed algorithm ADMM-LQP,  $r_{i}>0$ plays the role of proximal parameter,  $d(\cdot,\cdot)$  is a specified LQP regularizer given by  Definition  \ref{pre-1}  and the dual stepsizes
\begin{equation}\label{setK}
  (\alpha, \tau) \in \C{K}: = \left\{  (\alpha, \tau) \ | \ 1>\alpha>-1,~   \alpha + \tau >0,  \  1 + \alpha + \tau  -\alpha \tau -\alpha^2 - \tau^2 >0 \right\}.
\end{equation}
\begin{flushleft}
\centering\fbox{
	\parbox{0.89\textwidth}{
\[
   \begin{array}{lll}
 \hspace*{-.1in}\textrm{\bf Initialize } (x_1^0,\cdots,x_p^0, y^0,\lambda^0) \textrm{ and choose }
 (\alpha, \tau) \in \C{K},\beta>0; \\
 \hspace*{-.1in} \textrm{\bf While } \textrm{stopping criteria is not satisfied }  \textrm{\bf do}\\
\hspace*{.1in}  \textrm{\bf For}\ i=1,2,\cdots,p,~  \textrm{\bf do}\\
\hspace*{.3in}  x_i^{k+1}\leftarrow\arg\min\limits_{x_i\in \mathbb{R}_{+}^{m_i}} \C{L}_{\beta}(x_1^k,\cdots,x_{i-1}^k,x_i,x_{i+1}^k\cdots,x_p^k,y^{k},\lambda^{k})+r_id(x_i,x_i^k); \\
\hspace*{.1in}   \textrm{\bf End for}  \\
\hspace*{.1in}   \lambda^{k+\frac{1}{2}}\leftarrow\lambda^k-\alpha\beta(A \m{x}^{k+1}
+By^k-b);\\
\hspace*{.1in}     y^{k+1}\leftarrow\arg\min\limits_{y\in \C{Y}} \C{L}_{\beta}(x_1^{k+1},\cdots,x_p^{k+1},y,\lambda^{k+\frac{1}{2}}); \\
\hspace*{.1in}  \lambda^{k+1}\leftarrow\lambda^{k+\frac{1}{2}}-\tau\beta(A \m{x}^{k+1}
+By^{k+1}-b);\\
\hspace*{.1in}     k \leftarrow k+1;\\
 \hspace*{-.1in} \textrm{\bf End while}
\end{array}
\]}}\end{flushleft}\medskip

 The rest of this paper is organized as follows. In Section \ref{sec2}, we summarize some fundamental preliminaries including the LQP regularization term and its related results, the first-order optimality condition of the problem whose objective function is a smooth function plus nonsmooth function, {as well as} the variational characterization for the saddle point of the problem and the iterates of the {proposed} algorithm.  In Section \ref{sec3}, we  analyze the   global convergence  {of    ADMM-LQP with   the worst-case $O(1/T)$ convergence rate} for a new average iterate $\m{w}_{T}:=\frac {1}{1+T}\sum_{k=\kappa}^{\kappa+T}\tilde{\m{w}}^k $, where   $ \kappa\geq0, T>0$ are integers. Section \ref{sec445p} discusses and analyzes an extended version of ADMM-LQP. Finally, we  conclude    the paper in Section \ref{condr}.

{\bf Notations }  Let $\mathbb{R}, \mathbb{R}^n,\mathbb{R}^{m\times n}$
be the sets of  real numbers,  $n$ dimensional real column vectors and
 $m\times n$ {dimensional} real matrices, respectively. {In particular,} $\mathbb{R}_{++}^{n}(\mathbb{R}_{+}^{n})$ denotes the set of $n$  {dimensional real} positive (nonnegative) vectors.   {The bold  $\m{I}$  denotes  the identity matrix and $\m{0}$   denotes zero matrix/vector with proper dimensions.} For any symmetric matrix $G$, {we define $\|x\|_G^2 := x \tr Gx$}
where the superscript $\tr$ denotes the transpose.  Note that $G$ could be indefinite with $x \tr G x < 0$ for some $x$.
If $G$ is positive definite, then {we call $\|x\|_G$ $G$-weighted} norm. {The notations $\|\cdot\|, \langle \cdot,\cdot\rangle$ and $\nabla_x f(x,y)$
  denote  the standard Euclidean norm, inner product, and      the partial derivative of the differentiable function $f(x,y)$  at $x$, respectively.}
\section{Preliminaries} \label{sec2}

 The forthcoming {Definition  \ref{pre-1} and Lemma  \ref{lem1-b}, which are very} useful for analyzing convergence properties
 of ADMM-LQP in the sequel   sections, can be respectively  found in the earlier work \cite{ATe99} and \cite{LiYua15}.
It is clear {from   Definition  \ref{pre-1} that $d(\cdot, z)$ is a closed proper convex function and $d(\cdot, \cdot)$ is nonnegative. $d(v,z)=0$   if and
only if  $v=z$. Moreover,
\[
\nabla_{v}d(v,z)=(v-z)+\mu(z-Z^2v^{-1}),
\]
where $Z=\operatorname{diag}(z_1,z_2,\cdots,z_n)\in\mathbb{R}^{n\times n}$, $v^{-1}\in\mathbb{R}^n$ denotes a vector whose $j$th element is $1/v_j$.}
\begin{definition}\label{pre-1}
Let $0<\mu<1$ be a  given  constant. For $z \in{\mathbb{R}_{++}^{n}}$,  define {LQP regularizer as}
\[
d(v,z):=\begin{cases}
 \sum\limits_{j=1}^{n} \left[\frac {1}{2}\left(v_j-z_j\right)^2+\mu\left(z_j^2\log\frac {z_j}{v_j}+v_jz_j-z_j^2\right)\right],&\operatorname{if }\; v\in\mathbb{R}_{++}^{n},\\
{+\infty},&\operatorname{otherwise}.
\end{cases}
\]
\end{definition}

\begin{lemma}\label{lem1-b}
Let $P:=\operatorname{diag}(p_1,p_2,\cdots,p_n)\in\mathbb{R}^{n\times n}$ be a positive definite diagonal matrix, $q:\mathbb{R}_{+}^n\rightarrow\mathbb{R}^n$ be a monotone mapping and $\theta:\mathbb{R}^n\rightarrow\mathbb{R}$. Let $\mu$ be a given  positive constant. For given $\overline{z},z\in\mathbb{R}_{++}^n$, we define $\overline{Z}:=\operatorname{diag}(\overline{z}_1,\overline{z}_2,\cdots,\overline{z}_n)$, $z^{-1}:=(1/z_1,1/z_2,\cdots,1/z_n)\tr $ and
$
\phi^{'}(\bar{z},z):=(z-\bar{z})+\mu(\bar{z}-\bar{Z}^2z^{-1}).
$
Then, the following variational inequality
\[
\theta(v)-\theta(z)+ \left\langle v-z,  q(z)+P\phi^{'}(\bar{z},z)\right\rangle\geq 0,~~\forall v\in\mathbb{R}_{+}^n,
\]
has a unique positive solution $z$. Moreover,  for  a positive solution $z\in\mathbb{R}_{++}^n$, we have
\[
\theta(v)-\theta(z)+\left\langle v-z, q(z) \right\rangle\geq(1+\mu)\left\langle\bar{z}-z, P(v-z)\right\rangle-\mu\left\|\bar{z}-z\right\|_{P}^2,\quad \forall v\in\mathbb{R}_{+}^n.
\]
\end{lemma}

The following   lemma   can be found in \cite{HeMaYu016} that is {widely used to characterize} the first-order optimality conditions of  the subproblems  in   ADMM-LQP.
\begin{lemma}  \label{opt-1}
Let $f:\mathbb{R}^m\longrightarrow \mathbb{R}$ and
$h:\mathbb{R}^m\longrightarrow \mathbb{R}$ be two convex functions  defined on a closed convex
set $\Omega\subset \mathbb{R}^m$ and $h$ is differentiable. Suppose that the solution set
$\Omega^* = \arg\min\limits_{x\in\Omega}\{f(x)+h(x) \}$ is nonempty. Then, we have
\[
x^* \in \Omega^* \ \mbox{ if and only if } \ x^*\in \Omega, \  f(x)-f(x^*)
+\left\langle x-x^*, \nabla h(x^*)\right\rangle\geq 0,~ \forall x\in\Omega.
\]
\end{lemma}

Now, a point  $\m{w}^*=(\m{x}^*,y^*,\lambda^*)\in\ \C{W}$ is called the saddle-point of
  (\ref{Sec1-prob})
if
\[
\C{L}(\m{x}^*,y^*,\lambda)
\leq\C{L}(\m{x}^*,y^*,\lambda^*)
\leq \C{L}(\m{x},y,\lambda^*),~~ \forall  (\m{x}, y, \lambda) \in \C{W},
\]
  which is equivalent to
\[
  \left\{\begin{array}{lll}
 \textrm{\bf For}\ i=1,2,\cdots,p,  \\
   x_i^*\in \mathbb{R}_{+}^{m_i},~  f_i(x_i)- f_i(x_i^*) + \left\langle x_i-x_i^*, -A_i\tr \lambda^*\right\rangle \geq  0, ~ \forall x_i\in\mathbb{R}_{+}^{m_i},\\
y^*\in\C{Y}^*, \quad g(y)- g(y^*) +\left\langle y-y^*, -B\tr \lambda^*\right\rangle\geq 0,\qquad  \forall y\in\C{Y},\\
\lambda^*\in \mathbb{R}^{n},\quad \left\langle\lambda-\lambda^*,   A_1x_1^*+A_2x_2^*+\cdots+A_px_p^*+By^*- b  \right\rangle\geq  0,~ \forall \lambda\in \mathbb{R}^{n}.
\end{array}\right.
\]
Rewriting these inequalities in a  compact Variational Inequality (VI) form, we have
\begin{equation}\label{Sec31-1}
  \quad \Theta(\m{w})- \Theta(\m{w}^*) +\left\langle\m{w} -\m{w}^*, \C{J}(\m{w}^*)\right\rangle \geq  0, \ \forall \m{w}\in \C{W},
\end{equation}
where
\begin{equation}\label{Sec31-1j} \Theta(\m{w})=\sum\limits_{i=1}^p f_i(x_i)+g(y), ~
\m{w}=\left(
\begin{array}{c}
   \m{ x} \\
    y \\
    \lambda
  \end{array}
  \right),~
\C{J}(\m{w})=\left(
             \begin{array}{c}
              -A\tr \lambda \\
               -B\tr \lambda \\
               A\m{x}+By-b \\
             \end{array}
           \right).
\end{equation}\smallskip

{\noindent{\bf Remark  1} \it
Let  $\C{W}^*$ be the solution set of (\ref{Sec31-1}). Then,
     ADMM-LQP can be also used to solve the following structured variational inequality problem:
\begin{equation}\label{Sec31-001}
\textrm{find }  \m{w}^*\in \C{W}^*    \textrm{ such that } \left\langle\m{w} -\m{w}^*, \hat{\C{J}}(\m{w}^*)\right\rangle \geq  0, \ \forall \m{w}\in \C{W},
\end{equation}
from which
\begin{equation}\label{Sec31-001j}
\hat{\C{J}}(\m{w})=\left(
             \begin{array}{c}
             \theta_1(\m{ x}) -A\tr \lambda \\
              \theta_2(y) -B\tr \lambda \\
               A\m{x}+By-b \\
             \end{array}
           \right),~ \textrm{and}~ \theta_1(\m{ x})\in \partial \sum\limits_{i=1}^p f_i(x_i),~ \theta_2(y)\in \partial g(y).
\end{equation}
Here, $\partial g$ denotes its   limiting-subdifferential.  If   $f_i$ and $g$ are   differentiable, then solving Problem (\ref{Sec1-prob}) is equivalent to  solving
(\ref{Sec31-001}) with $\theta_1(\m{ x})=  \sum\limits_{i=1}^p \nabla f_i(x_i), \theta_2(y)= \nabla g(y)$.
Besides, the inequality (\ref{Sec31-1}) can be  also   rewritten as
\begin{equation}\label{Sec31-1j1}
 \textrm{VI}(\Theta,\mathcal{J},\C{W}): ~ \Theta(\m{w})- \Theta(\m{w}^*) +\left\langle\m{w} -\m{w}^*, \C{J}(\m{w})\right\rangle \geq  0, \ \forall \m{w}\in \C{W}
\end{equation}
because the affine mapping $\C{J}$ is skew-symmetric and   satisfies
\begin{equation}\label{Sec2-6}
\left\langle \m{w}-\bar{\m{w}}, \C{J}(\m{w})-\C{J}(\bar{\m{w}})\right\rangle\equiv0 ~~
\forall \m{w},\bar{\m{w}}\in\C{W}.
\end{equation}
Due to the  assumption that   the solution set of (\ref{Sec1-prob}) is nonempty, the solution set $\C{W}^*$   is also nonempty and convex. Moreover, it can be expressed as (see \cite{hy512}):
\[
\C{W}^*=\bigcap \limits_{\m{w}\in\C{W}}\left\{ \widehat{\m{w}} \in\C{W} ~| ~
\Theta(\m{w})- \Theta(\widehat{\m{w}}) +\left\langle\m{w} -\widehat{\m{w}}, \C{J}(\m{w})\right\rangle \geq  0 \right\}.
\]
 }

Note that  if the iterates generated by   ADMM-LQP  {satisfy}    $\textrm{VI}(\Theta,\mathcal{J},\C{W})$   with an extra term converging to {a fixed value as $k$ goes to infinity,
then global} convergence of our ADMM-LQP could be demonstrated theoretically. To verify this conjecture,
  we  introduce  the following auxiliary notations to simplify analysis:
\begin{equation}\label{Sec2-8}
                     \tilde{\m{x}}^k=\left(
              \begin{array}{c}
                \tilde{x}_{1}^{k}\\
                \tilde{x}_{2}^{k} \\
                \vdots\\
                \tilde{x}_{p}^{k} \\
              \end{array}
            \right)=\left(
                      \begin{array}{c}
                        x_{1}^{k+1} \\
                        x_{2}^{k+1} \\
                        \vdots \\
                        x_{p}^{k+1} \\
                      \end{array}
                    \right),~ \tilde{y}^k=y^{k+1},~
                    \tilde{\m{w}}^k=\left(
              \begin{array}{c}
                \tilde{\m{x}}^k \\
                \tilde{y}^k  \\
                \tilde{\lambda}^k  \\
              \end{array}
            \right)=\left(
                      \begin{array}{c}
                        \m{x}^{k+1} \\
                        y^{k+1}\\
                        \tilde{\lambda}^k  \\
                      \end{array}
                    \right),
\end{equation}
where
\begin{equation}\label{Sec2-9}
\tilde{\lambda}^k=\lambda^k-\beta(A\m{x}^{k+1}+By^{k}-b).
\end{equation}

\begin{lemma} (Prediction step)\label{le3t}
The iterates generated by ADMM-LQP   satisfy
\begin{equation}\label{Sec2-11}
\Theta(\m{w})-\Theta(\tilde{\m{w}}^k)+
\left\langle \m{w}-\tilde{\m{w}}^k,  \C{J}(\tilde{\m{w}}^k)+
Q(\tilde{\m{w}}^k-\m{w}^k)\right\rangle+\left\|\m{w}^k-\tilde{\m{w}}^k\right\|_{N}^2\geq0
\end{equation}
for any $\m{w}\in\C{W}$ where
\[
N=\begin{bmatrix}
      N_1 & \m{0 }\\
      \m{0} & \m{0 }\\
    \end{bmatrix}
\]
with {$N_1=\mu\operatorname{diag}( r_{1}\m{I},  r_{2}\m{I},\cdots, r_{p}\m{I})$},
and
\begin{equation}\label{Sec2-13}
Q=\begin{bmatrix}
      Q_1 & \m{0 }\\
      \m{0} & Q_2\\
    \end{bmatrix}
\end{equation}
{with}
\begin{equation}\label{Sec2-14}
Q_1=\begin{bmatrix}
       (1+\mu)r_{1}\m{I} & -\beta A_{1}\tr  A_{2} & \cdots& -\beta A_{1}\tr  A_{p} \\
        -\beta A_{2}\tr  A_{1} & (1+\mu)r_{2}\m{I} &  \cdots & -\beta A_{2}\tr  A_{p} \\
        \vdots & \vdots & \ddots & \vdots \\
        -\beta A_{p}\tr  A_{1} & -\beta A_{p}\tr  A_{2} & \cdots & (1+\mu)r_{p}\m{I} \\
      \end{bmatrix},~
Q_2=\begin{bmatrix}
       \beta B\tr B &  -\alpha B\tr  \\
        -B &  \frac 1\beta \m{I}\\
      \end{bmatrix}.
\end{equation}
\end{lemma}
\textbf{Proof}.   Using   {Lemma \ref{opt-1} and previous notations}, the first-order optimality condition of
 $x_{i}$-subproblem is
 \begin{eqnarray} \label{xi-opt}
&&\tilde{x}_{i}^k\in \mathbb{R}_{+}^{m_i},\quad f_{i}(x_{i})-f_{i}(\tilde{x}_{i}^k)+ \left\langle x_{i}-\tilde{x}_{i}^{k},
-A_{i}\tr\lambda^k+\right.\nonumber\\
&&\left.\beta A_{i}\tr\left(A_{i}\tilde{x}_{i}^{k}
+\sum_{r=1,r\neq i}^{p}A_{r}x_{r}^k+ B y^{k}-b\right) +r_{i}\nabla_{x_{i}^{k}} d(x_{i}^{k},\tilde{x}_{i}^k)\right\rangle\geq0
\end{eqnarray}
for any $x_i \in \mathbb{R}_{+}^{m_i}.$ Applying Lemma \ref{lem1-b} with
 $
 \theta(v)=f_i(x_i), v=x_i, z= \tilde{x}_{i}^k, \bar{z}=x_i^k, P= r_i \m{I}, q(z)=-A_{i}\tr\lambda^k+
 \beta A_{i}\tr\left(A_{i}\tilde{x}_{i}^{k}
+\sum_{r=1,r\neq i}^{p}A_{r}x_{r}^k+ B y^{k}-b\right)
 $
to (\ref{xi-opt}), we obtain
 \begin{eqnarray} \label{Sec2-16}
&&f_{i}(x_{i})-f_{i}(\tilde{x}_{i}^k)+ \left\langle x_{i}-\tilde{x}_{i}^{k},
-A_{i}\tr\tilde{\lambda}^k+ \beta A_{i}\tr \sum\limits_{r=1,r\neq i}^{p}A_{r}(x_{r}^k-\tilde{x}_{r}^k)\right\rangle\nonumber\\
&&\geq(1+\mu)r_{i}\left\langle x_{i}-\tilde{x}_{i}^{k},x_{i}^k-\tilde{x}_{i}^k\right\rangle
-\mu r_{i}\left\|x_{i}^k-\tilde{x}_{i}^k\right\|^2.
\end{eqnarray}
Meanwhile, the first-order optimality condition of $y$-subproblems  is
\begin{equation}\label{y-op-j1}
{\tilde{y}^k}\in \C{Y},\quad g(y)-g(\tilde{y}^k)+\left\langle y-\tilde{y}^{k},
-B\tr\lambda^{k+\frac 12}+\beta B\tr\left(A\tilde{\m{x}}^k  +B\tilde{y}^{k}-b\right)\right\rangle\geq0
\end{equation}
for any $ y\in \C{Y}$. By the way of generating $\lambda^{k+\frac 12}$   and    (\ref{Sec2-9}), it holds
\begin{equation}\label{Sec2-17}
\lambda^{k+\frac 12}=\lambda^{k}-\alpha(\lambda^{k}-\tilde{\lambda}^{k})=\tilde{\lambda}^{k}+(1-\alpha)
(\lambda^{k}-\tilde{\lambda}^{k}).
\end{equation}
{Plugging it into   (\ref{y-op-j1}) and using} an equivalent reformulation of (\ref{Sec2-9}), i.e.,
\begin{equation}\label{equ-refor}
\beta(A\tilde{\m{x}}^k +B\tilde{y}^{k}-b)=\lambda^{k}-\tilde{\lambda}^{k}+\beta B(\tilde{y}^k-y^k),
\end{equation}
{are to achieve}
\begin{equation}\label{Sec2-18}
g(y)-g(\tilde{y}^k)+\left\langle y-\tilde{y}^{k},
-B\tr\tilde{\lambda}^{k}+\alpha B\tr (\lambda^{k}-\tilde{\lambda}^{k})
+ \beta B\tr B(\tilde{y}^{k}-y^{k})\right\rangle\geq0.
\end{equation}
Notice that, the equation (\ref{equ-refor}) is equivalent to
\begin{equation}\label{Sec2-19}
\left\langle \lambda-\tilde{\lambda}^{k},  A\tilde{\m{x}}^k +B\tilde{y}^{k}-b
+\frac{1}{\beta}(\tilde{\lambda}^{k}-\lambda^{k})-  B(\tilde{y}^{k}-y^{k})
\right\rangle\geq0
\end{equation}
for any $\lambda\in \mathbb{R}^n$. Combining these inequalities (\ref{Sec2-16}), (\ref{Sec2-18})-(\ref{Sec2-19}) and the  definitions of matrices $N$ and $Q$, the   proof is completed. $\blacksquare$

Lemma \ref{le3t} shows that    the iterates generated by ADMM-LQP   can be characterized {as} a variational inequality with the aid of  an immediate variable $\tilde{\m{w}}$. We call  $\tilde{\m{w}}$ the predicting variable and $ \m{w}^{k+1} $  the correcting variable. Moreover, they satisfy the following relationship:
\begin{lemma} (Correction step) \label{M-strc}
The iterates $\tilde{\m{w}}$ {defined} by (\ref{Sec2-8}) and
  $ \m{w}^{k+1} $   generated by ADMM-LQP   satisfy
\begin{equation}\label{Sec2-20}
\m{w}^{k+1}=\m{w}^k-M(\m{w}^k-\tilde{\m{w}}^k),
\end{equation}
where
\[
M= \begin{bmatrix}
      \m{I}& \m{0} & \m{0}\\
     \m{0} & \m{I} & \m{0}\\
    \m{0} & -\tau\beta B  & (\alpha+\tau)\m{I}
    \end{bmatrix}.
\]
\end{lemma}
\textbf{Proof}.
By the   {updates of $\lambda^{k+1}$  and  $\lambda^{k+\frac{1}{2}}$}, it can be deduced that
\[
\begin{split}
 \lambda^{k+1}&=\lambda^{k+\frac{1}{2}}-\tau\beta(A\m{x}^{k+1}+By^{k+1}-b)\\
&=\lambda^{k+\frac{1}{2}}-\tau\beta(A\m{x}^{k+1}+By^{k}-b)+\tau\beta B(y^k-y^{k+1})\\
&=\lambda^{k+\frac{1}{2}}-\tau(\lambda^k-\tilde{\lambda}^k)+\tau\beta B(y^k-y^{k+1})\\
&=\lambda^k-\alpha(\lambda^k-\tilde{\lambda}^k)-\tau(\lambda^k-\tilde{\lambda}^k)+\tau\beta B(y^k-y^{k+1})\\
&=\lambda^k-[(\alpha+\tau)(\lambda^k-\tilde{\lambda}^k)-\tau\beta B(y^k-\tilde{y}^k)].
\end{split}
\]
The above equality together with (\ref{Sec2-8}) immediately implies (\ref{Sec2-20}). $\blacksquare$
\section{Convergence analysis of ADMM-LQP}\label{sec3}
\subsection{\textbf{Basic properties of $\{\m{w}^k-\m{w}^*\}$}}

In the following, we begin to analyze some properties of the sequences $\{\m{w}^k-\m{w}^*\}$ under  a  special  $H$-weighted norm.
Followed by (\ref{Sec2-6}), (\ref{Sec2-11}) can be rewritten as
\begin{equation}\label{Sec3-23}
\Theta(\m{w})-\Theta(\tilde{\m{w}}^k)+\langle\m{w}-\tilde{\m{w}}^k, \C{J}(\m{w})\rangle
\geq\left\langle \m{w}-\tilde{\m{w}}^k, Q(\m{w}^k-\tilde{\m{w}}^k)\right\rangle-\left\|\m{w}^k-\tilde{\m{w}}^k\right\|_{N}^2.
\end{equation}
For any $\alpha+\tau>0$, the matrix $M$ involved in Lemma \ref{M-strc} is   nonsingular. So, by (\ref{Sec2-20}) and a direct calculation, the first term on the right-hand side of (\ref{Sec3-23}) {becomes}
\begin{equation}\label{Sec3-24}
\left\langle \m{w}-\tilde{\m{w}}^k, Q(\m{w}^k-\tilde{\m{w}}^k)\right\rangle=
\left\langle \m{w}-\tilde{\m{w}}^k, H(\m{w}^k-\m{w}^{k+1})\right\rangle
\end{equation}
with
\begin{equation}\label{Sec3-25}
H=QM^{-1}=\begin{bmatrix}
              Q_{1} & \m{0} \\
              \m{0} & H_{2} \\
            \end{bmatrix},
\end{equation}
where $Q_{1}$ is {given by} (\ref{Sec2-14}) and
\[
H_{2}=\begin{bmatrix}
          (1-\frac {\alpha\tau}{\alpha+\tau})\beta B\tr B & -\frac \alpha{\alpha+\tau}B\tr  \\
          -\frac \alpha{\alpha+\tau}B  & \frac {1}{(\alpha+\tau)\beta} \m{I} \\
        \end{bmatrix}.
\]

Next, we provide a sufficient condition to ensure the positive definiteness of  $H$.
\begin{lemma} \label{le5-j}
For any $0<\mu<1$ and {$\gamma>\frac {p-1}{1+\mu}$,
the} matrix $H$ defined in (\ref{Sec3-25}) is symmetric positive definite if
\begin{equation}\label{Sec3-27}
{r_{i}\geq\gamma \beta \left\|A_{i}\tr A_{i}\right\|,i=1,2,\cdots,p,~~ \alpha<1 ~\textrm{and}  ~\alpha+\tau>0}
\end{equation}
\end{lemma}
\textbf{Proof}.
It follows from the known conditions that
\[
H
            \succeq
           \begin{bmatrix}
              \tilde{Q}_{1} & \m{0} \\
              \m{0} & H_{2} \\
            \end{bmatrix} :=\tilde{H},
\]
where
\[
\begin{split}
\tilde{Q}_1=
        \beta\begin{bmatrix}
      (1+\mu)\gamma A_{1}\tr  A_{1} & - A_{1}\tr  A_{2} & \cdots & -  A_{1}\tr  A_{p} \\
        -  A_{2}\tr  A_{1} & (1+\mu)\gamma  A_{2}\tr  A_{2} &  \cdots & -  A_{2}\tr  A_{p} \\
        \vdots & \vdots & \ddots & \vdots \\
        -  A_{p}\tr  A_{1} & -  A_{p}\tr  A_{2} & \cdots & (1+\mu)\gamma  A_{p}\tr  A_{p} \\
      \end{bmatrix}.
\end{split}
\]
So,   $H$ is   positive definite  if  $\tilde{H}$ is positive definite. By the block structure  of $\tilde{H}$, we only  need to show that  both $\tilde{Q}_{1}$  and  $H_{2}$   are positive definite. Note that $\tilde{Q}_{1}$ can be decomposed as
\begin{equation}\label{Sec3-30}
\tilde{Q}_1=\beta{\begin{bmatrix}
                     A_{1} &  &  &  \\
                      & A_{2} &  &  \\
                      &  & \ddots &  \\
                     &  &  & A_{p} \\
                   \end{bmatrix}}\tr \tilde{Q}_{1,0}\begin{bmatrix}
                     A_{1} &  &  &  \\
                      & A_{2} &  &  \\
                      &  & \ddots &  \\
                     &  &  & A_{p} \\
                      \end{bmatrix}
\end{equation}
where
\[
\tilde{Q}_{1,0}={\begin{bmatrix}
                     \gamma(1+\mu)\m{I} & -\m{I}  & \cdots &-\m{I}  \\
                      -\m{I} & \gamma(1+\mu)\m{I} &\cdots  &-\m{I}  \\
                     \vdots &\vdots  & \ddots&  \vdots\\
                     -\m{I}& -\m{I} &\cdots  & \gamma(1+\mu)\m{I} \\
                      \end{bmatrix}} \textcolor{red}{.}
\]
For any $\gamma>\frac {p-1}{1+\mu}$, then $\tilde{Q}_{1,0}$ is  strictly diagonally dominant and thus   positive definite. So,     by  (\ref{Sec3-30})  $\tilde{Q}_1$ is positive definite {under the conditions that} $\gamma>\frac {p-1}{1+\mu}$ and all $A_{i}(i=1,\ldots,p)$ have full column ranks. In addition, we have
\[
H_{2}={\begin{bmatrix}
                  \beta^{\frac 12}B &  \\
                   & \beta^{-\frac 12}\m{I} \\
                \end{bmatrix}}\tr H_{2,0}\begin{bmatrix}
                  \beta^{\frac 12}B &  \\
                   & \beta^{-\frac 12}\m{I} \\
                \end{bmatrix}
\]
where
\[
H_{2,0}=\frac{1}{\alpha+\tau}\begin{bmatrix}
          (\alpha+\tau- \alpha\tau )\m{I} & - \alpha \m{I}  \\
          - \alpha \m{I} &  \m{I} \\
        \end{bmatrix}.
\]
$H_{2,0}$ is positive definite if  $(\alpha,\tau)$ satisfy (\ref{Sec3-27}).
 So, $H_{2}$ is positive definite if  the last two conditions of (\ref{Sec3-27}) hold   and the matrix $B$ has full column rank.

 Summarizing the above discussions, if the conditions in  (\ref{Sec3-27}) hold, then the matrix $H$ is positive definite. This completes the whole proof. $\blacksquare$
\begin{theorem}\label{t1-bjc}
The iterates generated by ADMM-LQP   satisfy
\begin{equation}\label{Sec3-34}
\begin{split}
&\Theta(\m{w})-\Theta(\tilde{\m{w}}^k)+\left\langle \m{w}-\tilde{\m{w}}^k, \C{J}(\m{w})\right\rangle\\
&\geq\frac {1}{2}\left\{\left\|\m{w}-\m{w}^{k+1}\right\|_{H}^2-\left\|\m{w}-\m{w}^k\right\|_{H}^2\right\}  +\frac {1}{2}\left\|\m{w}^k-\tilde{\m{w}}^k\right\|_{G}^2,~~\forall \m{w}\in\C{W},
\end{split}
\end{equation}
where
\begin{equation}\label{Sec3-35}
G=Q\tr +Q-M\tr HM-2N.
\end{equation}
\end{theorem}
\textbf{Proof}.
Taking
$
(a,b,c,d)=(\m{w},\tilde{\m{w}}^k,\m{w}^k,\m{w}^{k+1})
$
in  the following identity
\[
(a-b)\tr H(c-d)=\frac {1}{2}\left\{\|a-d\|_{H}^2-\|a-c\|_{H}^2\right\}+\frac {1}{2}\left\{\|c-b\|_{H}^2-\|d-b\|_{H}^2\right\}
\]
gives
\begin{equation}\label{Sec3-36}
\begin{split}
(\m{w}-\tilde{\m{w}}^k)\tr H(\m{w}^k-\m{w}^{k+1})&=\frac {1}{2}\left\{\left\|\m{w}-\m{w}^{k+1}\right\|_{H}^2-\left\|\m{w}-\m{w}^k\right\|_{H}^2\right\}\\
&\quad+\frac {1}{2}\left\{\left\|\m{w}^k-\tilde{\m{w}}^k\right\|_{H}^2-\left\|\m{w}^{k+1}-\tilde{\m{w}}^k\right\|_{H}^2\right\}.
\end{split}
\end{equation}
By a simple deduction, we have
\begin{equation}\label{Sec3-37}
\begin{split}
&\frac {1}{2}\left\{\left\|\m{w}^k-\tilde{\m{w}}^k\right\|_{H}^2-\left\|\m{w}^{k+1}-\tilde{\m{w}}^k\right\|_{H}^2\right\}
-\left\|\m{w}^k-\tilde{\m{w}}^k\right\|_{N}^2\\
&=\frac {1}{2}\left\{\left\|\m{w}^k-\tilde{\m{w}}^k\right\|_{H}^2-\left\|\m{w}^{k+1}-\m{w}^k
+\m{w}^k-\tilde{\m{w}}^k\right\|_{H}^2\right\}-\left\|\m{w}^k-\tilde{\m{w}}^k\right\|_{N}^2\\
&=\frac {1}{2}\left\{\left\|\m{w}^k-\tilde{\m{w}}^k\right\|_{H}^2-\left\|\m{w}^k-\tilde{\m{w}}^k
-M(\m{w}^k-\tilde{\m{w}}^k)\right\|_{H}^2\right\}-\left\|\m{w}^k-\tilde{\m{w}}^k\right\|_{N}^2\\
&=\frac {1}{2}\left\langle \m{w}^k-\tilde{\m{w}}^k, \left(HM+(HM)\tr-M\tr HM-2N\right)(\m{w}^k-\tilde{\m{w}}^k)\right\rangle\\
&=\frac {1}{2}\left\|\m{w}^k-\tilde{\m{w}}^k\right\|_{Q\tr +Q-M\tr HM-2N}^2,
\end{split}
\end{equation}
where the second and the last \textcolor{blue}{equalities} follow from (\ref{Sec2-20})   and (\ref{Sec3-25}), respectively. Finally, (\ref{Sec3-34}) follows from (\ref{Sec3-23})-(\ref{Sec3-24}), (\ref{Sec3-36})-(\ref{Sec3-37}) and the definition of $G$ in (\ref{Sec3-35}). $\blacksquare$

In fact, setting $\m{w}=\m{w}^*\in\C{W}^*$ in (\ref{Sec3-34})  we could have
\[
\begin{split}
&\frac {1}{2}\left\{\left\|\m{w}^*-\m{w}^k\right\|_{H}^2-\left\|\m{w}^*-\m{w}^{k+1}\right\|_{H}^2-
\left\|\m{w}^k-\tilde{\m{w}}^k\right\|_{G}^2\right\}\\
& \geq
\Theta(\tilde{\m{w}}^k)-\Theta(\m{w}^*)+\left\langle \tilde{\m{w}}^k-\m{w}^*, \C{J}(\m{w}^*)\right\rangle\geq 0,
\end{split}
\]
from which  the following theorem holds immediately.

\begin{theorem} \label{contr}
The iterates generated by ADMM-LQP satisfy
\begin{equation}\label{Sec3-38}
\left\|\m{w}^{k+1}-\m{w}^*\right\|_{H}^2\leq\left\|\m{w}^k-\m{w}^*\right\|_{H}^2
-\left\|\m{w}^k-\tilde{\m{w}}^k\right\|_{G}^2,~~ \forall \m{w}^*\in\C{W}^*.
\end{equation}
\end{theorem}

Observing that if the   matrix $G$ in Theorem \ref{contr} is positive definite, then similar convergence results to \cite{Bazl018,HeYuanLQP06,hy512,HeXuYuan2016b} {hold}  for  our ADMM-LQP. However, the   matrix $G$ is not necessarily  positive definite for any $ (\alpha, \tau)$ satisfying (\ref{setK}). Therefore, it is full of necessity to estimate the lower
bound of $\left\|\m{w}^k-\tilde{\m{w}}^k\right\|_{G}^2$ for further   analysis.

\subsection{\textbf{Lower bound of $\left\|\m{w}^k-\tilde{\m{w}}^k\right\|_{G}^2$}}\label{se32b}
In this subsection, we focus on estimating the lower bound of $\left\|\m{w}^k-\tilde{\m{w}}^k\right\|_{G}^2$ by combining the structure of $G$ as well as the first-order optimality {condition} of $y$-subproblem.
 By {a} direct calculation  we have
\begin{equation}\label{Sec3-35j}
G= \begin{bmatrix}
                               G_{1} & \m{0} \\
                               \m{0} & G_{2} \\
                             \end{bmatrix},
\end{equation}
where
\[
G_{1}=\begin{bmatrix}
          (1-\mu)r_{1}\m{I} & -\beta A_{1}\tr A_{2} & \cdots & -\beta A_{1}\tr A_{p} \\
          -\beta A_{2}\tr A_{1} & (1-\mu)r_{2}\m{I} & \cdots & -\beta A_{2}\tr A_{p} \\
          \vdots & \vdots & \ddots & \vdots \\
          -\beta A_{p}\tr A_{1} & -\beta A_{p}\tr A_{2} & \cdots & (1-\mu)r_{p}\m{I}  \\
        \end{bmatrix},
G_{2}=\begin{bmatrix}
          (1-\tau)\beta B\tr B  & (\tau-1)B\tr  \\
          (\tau-1)B & \frac {2-\alpha-\tau}{\beta}\m{I} \\
        \end{bmatrix}.
\]
In the following, we  define a set of auxiliary variables
 \begin{equation}\label{Sec2-8j0-1}
\m{E}^{k+1}=A \m{x}^{k+1}
+By^{k+1}-b,~~  \m{E}_{\m{x}}^{k+1}= \m{x}^{k+1}- \m{x}^k,~~\textrm{and}~~ \m{E}_y^{k+1}=y^{k+1}- y^k.
 \end{equation}

\begin{lemma}\label{le6-j}
  {For any $\gamma>\frac {p-1}{1-\mu}$, if $r_{i}\geq\gamma\beta \left\|A_{i}\tr A_{i}\right\|$ for $ i=1,2,\cdots,p$}, then    there exists a constant $\xi_{1}>0$ such that
\begin{equation}\label{Sec3-39}
\begin{split}
\left\|\m{w}^k-\tilde{\m{w}}^k\right\|_{G}^2\geq&\xi_{1}\sum\limits_{i=1}^p \left\|A_{i}\m{E}_{x_i}^{k+1}\right\|^2+\beta(1-\alpha)\left\|B\m{E}_y^{k+1}\right\|^2\\&
+\beta(2-\alpha-\tau)\left\|\m{E}^{k+1}\right\|^2
-2\beta(1-\alpha)\left\langle \m{E}^{k+1}, B\m{E}_y^{k+1}\right\rangle.
\end{split}
\end{equation}
\end{lemma}
\textbf{Proof}.
 { Since $r_{i}\geq\gamma\beta \left\|A_{i}\tr A_{i}\right\|$}, we deduce
\[
\begin{split}
G_{1}&=\begin{bmatrix}
          (1-\mu)r_{1}\m{I} & -\beta A_{1}\tr A_{2} & \cdots & -\beta A_{1}\tr A_{p} \\
          -\beta A_{2}\tr A_{1} & (1-\mu)r_{2}\m{I} & \cdots & -\beta A_{2}\tr A_{p} \\
          \vdots & \vdots & \ddots & \vdots \\
          -\beta A_{p}\tr A_{1} & -\beta A_{p}\tr A_{2} & \cdots & (1-\mu)r_{p}\m{I}  \\
        \end{bmatrix}\\&
        \succeq
        {\beta\begin{bmatrix}
          (1-\mu)\gamma A_{1}\tr A_{1} & -  A_{1}\tr A_{2} & \cdots & -  A_{1}\tr A_{p} \\
          -  A_{2}\tr A_{1} & (1-\mu)\gamma  A_{2}\tr A_{2} & \cdots & -  A_{2}\tr A_{p} \\
          \vdots & \vdots & \ddots & \vdots \\
          -  A_{p}\tr A_{1} & -  A_{p}\tr A_{2} & \cdots & (1-\mu)\gamma  A_{p}\tr A_{p}  \\
        \end{bmatrix}}=\tilde{G}_{1}.
\end{split}
\]
So,  it is clear that
\[
G \succeq
\begin{bmatrix}
              \tilde{G}_{1} & \m{0}\\
             \m{0 } & G_{2} \\
            \end{bmatrix}:=\tilde{G}\quad \Longrightarrow \quad
\left\|\m{w}^k-\tilde{\m{w}}^k\right\|_{G}^2\geq\left\|\m{w}^k-\tilde{\m{w}}^k\right\|_{\tilde{G}}^2.
\]
Then, by the structure of  $\tilde{G}$ the following equality holds:
\[
\begin{split}
\left\|\m{w}^k-\tilde{\m{w}}^k\right\|_{\tilde{G}}^2=
&\beta{\begin{Vmatrix}\begin{pmatrix}
                                                           A_{1}(x_{1}^k-x_{1}^{k+1}) \\
                                                           A_{2}(x_{2}^k-x_{2}^{k+1}) \\
                                                           \vdots\\
                                                           A_{p}(x_{p}^k-x_{p}^{k+1}) \\
                                                         \end{pmatrix}\end{Vmatrix}}
                                                       _{\bar{G}_{1,0}}^2
+\beta(1-\tau)\left\|B\m{E}_y^{k+1}\right\|^2
\\&+2(1-\tau)\left\langle\lambda^k-\tilde{\lambda}^k, B\m{E}_y^{k+1}\right\rangle+\frac {2-\alpha-\tau}{\beta}\left\|\lambda^k-\tilde{\lambda}^k\right\|^2,
  \end{split}
 \]
where
\[
\bar{G}_{1,0}={\begin{bmatrix}
                  \gamma(1-\mu)\m{I} & -\m{I} & \cdots & -\m{I} \\
                  -\m{I} & \gamma(1-\mu)\m{I} & \cdots & -\m{I} \\
                  \vdots & \vdots & \ddots & \vdots \\
                  -\m{I} & -\m{I} & \cdots & \gamma(1-\mu)\m{I} \\
                \end{bmatrix}}.
\]
Since {$0<\mu<1$, $\gamma>\frac {p-1}{1-\mu}$} and $\bar{G}_{1,0}$ is positive definite,   there exists a    $\xi_{1}>0$ such that
\[
\beta{\begin{Vmatrix}\begin{pmatrix}
                                                           A_{1}(x_{1}^k-x_{1}^{k+1}) \\
                                                           A_{2}(x_{2}^k-x_{2}^{k+1}) \\
                                                           \vdots\\
                                                           A_{p}(x_{p}^k-x_{p}^{k+1}) \\
                                                         \end{pmatrix}\end{Vmatrix}}
                                                       _{\bar{G}_{1,0}}^2
\geq\xi_{1}\sum\limits_{i=1}^p \left\|A_{i}\m{E}_{x_i}^{k+1}\right\|^2.
\]
The definition of $\tilde{\lambda}^k$ in (\ref{Sec2-9}) shows
$
\lambda^k-\tilde{\lambda}^k=\beta\left[ \m{E}^{k+1}-B \m{E}_y^{k+1}\right].
$
So, by   (\ref{Sec2-8j0-1}), we  see  that
\[
\begin{split}
 2(1-\tau)\left\langle\lambda^k-\tilde{\lambda}^k, \textcolor{blue}{B \m{E}_y^{k+1}}\right\rangle
=2\beta(1-\tau)\left\langle  \m{E}^{k+1}, B \m{E}_y^{k+1}\right\rangle
+2\beta(\tau-1)\left\|B \m{E}_y^{k+1}\right\|^2,
\end{split}
\]
and
\[
\begin{split}
& \frac {2-\alpha-\tau}{\beta}\left\|\lambda^k-\tilde{\lambda}^k\right\|^2=\beta(2-\alpha-\tau)\left\|\m{E}^{k+1}\right\|^2\\
&
+
\beta(2-\alpha-\tau)\left\|B\m{E}_y^{k+1}\right\|^2
  -2\beta(2-\alpha-\tau)\left\langle \m{E}^{k+1}, B\m{E}_y^{k+1}\right\rangle.
\end{split}
\]
The above discussions {illustrate that}
\[
\begin{split}
\left\|\m{w}^k-\tilde{\m{w}}^k\right\|_{\tilde{G}}^2\geq
&\xi_{1}\sum\limits_{i=1}^p \left\|A_{i}\m{E}_{x_{i}}^{k+1}\right\|^2
+\beta(1-\alpha)\left\|B\m{E}_y^{k+1}\right\|^2\\
&
+\beta(2-\alpha-\tau)\|\m{E}^{k+1}\|^2
-2\beta(1-\alpha)\left\langle \m{E}^{k+1}, B\m{E}_y^{k+1}\right\rangle.
\end{split}
\]
This completes the proof. $\blacksquare$

\begin{lemma}
Suppose $\alpha>-1$. Then the iterates   generated by  ADMM-LQP satisfy
\begin{equation}\label{Sec3-44}
\begin{split}
\left\langle \m{E}^{k+1}, -B\m{E}_y^{k+1}\right\rangle
\geq\frac {\tau-1}{1+\alpha}\left\langle \m{E}^k, B\m{E}_y^{k+1}\right\rangle
-\frac {\alpha}{1+\alpha}\left\|B\m{E}_y^{k+1}\right\|^2.
\end{split}
\end{equation}
\end{lemma}
\textbf{Proof}.
The optimality {condition  of $y$-subproblem is}
\[y^{k+1}\in \C{Y},\quad
g(y)-g(y^{k+1})+\left\langle y-y^{k+1},
-B\tr \lambda^{k+\frac {1}{2}}+\beta B\tr \m{E}^{k+1}\right\rangle\geq0
\]
for any $y\in \C{Y}$. Setting $y=y^k$ in the above inequality gives
\begin{equation}\label{Sec3-45}
g(y^k)-g(y^{k+1})-\left\langle \m{E}_y^{k+1},
-B\tr \lambda^{k+\frac 12}+\beta B\tr \m{E}^{k+1}\right\rangle\geq0.
\end{equation}
Similarly,   the optimality condition of {$y$-subproblem at the previous iteration is}
\[
g(y)-g(y^k)+\left\langle y-y^k,
-B\tr \lambda^{k-\frac {1}{2}}+\beta B\tr \m{E}^k\right\rangle\geq0.
\]
Letting $y=y^{k+1}$ in the above inequality, we have
\begin{equation}\label{Sec3-46}
g(y^{k+1})-g(y^k)+\left\langle \m{E}_y^{k+1},
-B\tr \lambda^{k-\frac {1}{2}}+\beta B\tr \m{E}^k\right\rangle\geq0.
\end{equation}
Then, adding (\ref{Sec3-45}) and (\ref{Sec3-46}) together is to obtain
\begin{equation}\label{Sec3-47}
\begin{split}
-\left\langle B\m{E}_y^{k+1},
   \lambda^{k-\frac 12} -\lambda^{k+\frac 12} +\beta ( \m{E}^{k+1}
-  \m{E}^k)\right\rangle\geq0.
\end{split}
\end{equation}
Notice by the updates of $\lambda^{k+\frac 12}$ and $\lambda^{k}$, i.e.,
$
\lambda^{k+\frac 12}=\lambda^k-\alpha\beta(A\m{x}^{k+1}+By^{k}-b)$
  and $
\lambda^{k}=\lambda^{k-\frac 12}-\tau\beta\m{E}^k
$,
that
\begin{eqnarray}\label{sjb0-0}
\lambda^{k-\frac 12}-\lambda^{k+\frac 12}&=&\tau\beta\m{E}^k
+ \alpha\beta\m{E}^{k+1}-\alpha\beta B\m{E}_y^{k+1}.
\end{eqnarray}
Substituting (\ref{sjb0-0}) into the \textcolor{blue}{left-hand} side  of (\ref{Sec3-47}), we get
\[
\begin{split}
& \left\langle-B\m{E}_y^{k+1},  (\alpha+1)   \m{E}^{k+1}
+ (\tau-1) \m{E}^k
 -\alpha   B\m{E}_y^{k+1}  \right\rangle\geq 0
\end{split}
\]
which is equivalent to (\ref{Sec3-44}).  $\blacksquare$

The following result provides  a concrete form of  the lower bound of $\left\|\m{w}^k-\tilde{\m{w}}^k\right\|_{G}^2$, which is the key technique to {analyze}  convergence properties of the proposed algorithm.
\begin{theorem} (Lower bound estimation)\label{Theor2}
Let the sequences $\{\m{w}^k\}$ be generated by ADMM-LQP and $\{\tilde{\m{w}}^k\}$  be defined by (\ref{Sec2-8}). {For any
\[
\gamma>\frac {p-1}{1-\mu} ~~ \textrm{and}~~(\alpha,\tau)\in \mathcal{K},
\]
where $0<\mu<1$ and $\C{K}$ is defined by (\ref{setK}),
  there exist   $\xi_{1},\xi_{2}>0, \xi_{3}\geq0$}  such that
\[
\begin{split}
\left\|\m{w}^k-\tilde{\m{w}}^k \right\|_{G}^2\geq&\xi_{1}\sum\limits_{i=1}^p \left\|A_{i}\m{E}_{x_{i}}^{k+1}\right\|^2
+\xi_{2}\left\|\m{E}^{k+1} \right\|^2
+\xi_{3}\left(\left\|\m{E}^{k+1}\right\|^2-\left\|\m{E}^k\right\|^2\right).
\end{split}
\]
\end{theorem}
\textbf{Proof}. First of all, it follows from (\ref{Sec3-39}) and (\ref{Sec3-44}) that
there exists a $\xi_{1}>0$ such that
\begin{equation}\label{Sec3-48}
\begin{split}
&\left\|\m{w}^k-\tilde{\m{w}}^k\right\|_{G}^2\geq \xi_{1}\sum\limits_{i=1}^p \left\|A_{i}\m{E}_{x_{i}}^{k+1}\right\|^2+\frac {\beta(1-\alpha)^2}{1+\alpha}\left\|B\m{E}_y^{k+1}\right\|^2\\
&+ \beta(2-\alpha-\tau)\left\|\m{E}^{k+1}\right\|^2
+\frac {2\beta(1-\alpha)(\tau-1)}{1+\alpha}\left\langle \m{E}^k, B\m{E}_y^{k+1}\right\rangle.
\end{split}
\end{equation}
By  the Cauchy-Schwartz inequality, the following holds
\[
\begin{split}
 2(1-\alpha)(\tau-1)\left\langle \m{E}^k, B\m{E}_y^{k+1}\right\rangle
\geq -(1-\tau)^2\left\|\m{E}^k\right\|^2-(1-\alpha)^2\left\|B\m{E}_y^{k+1}\right\|^2.
\end{split}
\]
 {Substituting it into   (\ref{Sec3-48})  gives}
\[
\begin{split}
&\left\|\m{w}^k-\tilde{\m{w}}^k \right\|_{G}^2\geq\xi_{1}\sum\limits_{i=1}^p \left\|A_{i}\m{E}_{x_{i}}^{k+1}\right\|^2\\
&
+\beta\left(2-\alpha-\tau-\frac {(1-\tau)^2}{1+\alpha}\right)\left\|\m{E}^{k+1}\right\|^2
+\frac {\beta(1-\tau)^2}{1+\alpha}\left(\left\|\m{E}^{k+1}\right\|^2-\left\|\m{E}^k\right\|^2\right).
\end{split}
\]
For any $(\alpha,\tau)\in \mathcal{K}$ and $\beta>0$, it holds that
\[
\xi_{2}:=\beta\left(2-\alpha-\tau-\frac {(1-\tau)^2}{1+\alpha}\right)>0\;  ~ \textrm{ and }~
\;\xi_{3}:=\frac {\beta(1-\tau)^2}{1+\alpha}\geq0,
\]
which confirms the conclusion. $\blacksquare$
 \subsection{\textbf{Global convergence and sublinear convergence rate}} \label{33glob}
In this part, we focus on analyzing the global convergence of the proposed   ADMM-LQP and establishing its sublinear convergence rate. The following corollary is obtained directly from the precious Theorems \ref{t1-bjc}, \ref{contr} and \ref{Theor2}, respectively.
\begin{corollary}
Suppose that conditions of Theorem \ref{Theor2} hold. Then,   there exist $\xi_{1},\xi_{2}>0$ and $\xi_{3}\geq0$ such that the iterates generated by ADMM-LQP satisfy
\begin{eqnarray}
\begin{split}
&\left\|\m{w}^{k+1}-\m{w}^*\right\|_{H}^2+\xi_{3}\left\|\m{E}^{k+1}\right\|^2 \label{Sec3-54}\\
& \leq\left\|\m{w}^k-\m{w}^*\right\|_{H}^2+\xi_{3}\left\|\m{E}^k\right\|^2
-\xi_{1}\sum\limits_{i=1}^p \left\|A_{i}\m{E}_{x_{i}}^{k+1}\right\|^2
-\xi_{2}\left\|\m{E}^{k+1}\right\|^2
\end{split}
\end{eqnarray}
and
\begin{eqnarray}\label{bj-ke1}
\begin{split}
&\Theta(\m{w})-\Theta(\tilde{\m{w}}^k)+\left\langle \m{w}-\tilde{\m{w}}^k, \C{J}(\m{w})\right\rangle
\\
&\geq\frac {1}{2}\left\{\left\|\m{w}-\m{w}^{k+1}\right\|_{H}^2+ \xi_{3}\left\|\m{E}^{k+1}\right\|^2\right\}
-\frac{1}{2}\left\{\left\|\m{w}-\m{w}^k\right\|_{H}^2+\xi_{3}\left\|\m{E}^k\right\|^2\right\}.
\end{split}
\end{eqnarray}
\end{corollary}

\begin{theorem} (Global convergence)\label{conver-th}
Suppose the conditions of Theorem \ref{Theor2} hold.  Then, \par
\noindent (i)
   $
\lim\limits_{k\rightarrow\infty}\sum\limits_{i=1}^p \left\|A_{i}(x_{i}^k-x_{i}^{k+1})\right\|^2=0 ~ \textrm{and}~
\lim\limits_{k\rightarrow\infty}\left\|A\m{x}^{k+1}+By^{k+1}-b\right\|^2=0;
$

\noindent (ii) any accumulation point of $\{\m{w}^k\}$ is a solution of $\textrm{VI}(\Theta,\mathcal{J},\C{W})$;\par

\noindent (iii) there exists a {point} $\m{w}^{\infty}\in\C{W}^*$ such that
$
\lim\limits_{k\rightarrow\infty}\m{w}^k=\m{w}^{\infty}.
$
\end{theorem}
\textbf{Proof}.  The statement (i)  can be proved by
summing  (\ref{Sec3-54}) over $k=0,1,\ldots,\infty$ directly:
\begin{equation}\label{Sec3-591ij}
\begin{split}
 \sum\limits_{k=0}^{\infty}\left(\xi_{1}\sum\limits_{i=1}^p \left\|A_{i}\m{E}_{x_{i}}^{k+1}\right\|^2
+\xi_{2}\left\|\m{E}^{k+1} \right\|^2\right)
\leq\left\|\m{w}^0-\m{w}^*\right\|_{H}^2+\xi_{3}\left\|\m{E}^0\right\|^2<\infty,
\end{split}
\end{equation}
since $\xi_{1},\xi_{2}>0$.

Now, we prove (ii). It follows from (i) ,    (\ref{Sec2-8j0-1}), and the full column rank assumption on all the matrices $A_{i}$   that
\begin{equation}\label{Sec3-591}
\lim_{k\rightarrow\infty} x_{i}^k-x_{i}^{k+1}=\m{0} \quad
\textrm{and}\quad
\lim_{k\rightarrow\infty}{A\m{x}^{k+1}+By^{k+1}-b}=\m{0}
\end{equation}
for all $i=1,\ldots,p$. By   the definitions of $\tilde{\m{x}}^k$ and $\tilde{\lambda}^k$ in (\ref{Sec2-8}) and (\ref{Sec2-9}), we have
\[
 \lambda ^k-\tilde{\lambda}^k= \beta A(\tilde{\m{x}}^k-\m{x}^k) +\beta(A\m{x}^k+By^k-b).
\]
Combining the above equation with (\ref{Sec3-591}), it holds
\[
\lim_{k\rightarrow\infty}{\lambda ^k-\tilde{\lambda}^k}=\m{0}.
\]
The above {discussions} together with   (\ref{equ-refor}) as well as the full column rank assumption on $B$ show
$
\lim\limits_{k\rightarrow\infty}{y^k-\tilde{y}^k}=\m{0}.
$
Hence, we have
\begin{equation}\label{Sec3-59bj1}
\lim_{k\rightarrow\infty}{\m{w}^k-\tilde{\m{w}}^k}=\m{0}.
\end{equation}
Suppose $\hat{\m{w}}=(\hat{\m{x}},\hat{y},\hat{\lambda})$ is an accumulation point of $\{\m{w}^k\}$, that is, there is a subsequence $\{\m{w}^{k_j}\}$ converging to $\hat{\m{w}}$ as $j$ goes to infinity. Then, we have
\[
\lim\limits_{j\rightarrow\infty}\m{w}^{k_j}  = \hat{\m{w}}=\lim\limits_{k\rightarrow\infty}\tilde{\m{w}}^{k}.
\]
So, for any fixed $\m{w}\in\C{W}$, taking $k=k_j$ in (\ref{Sec2-11}) and letting $j\rightarrow\infty$, we have \[
\Theta(\m{w})-\Theta(\hat{\m{w}})+
\left\langle \m{w}- \hat{\m{w}},  \C{J}(\hat{\m{w}}) \right\rangle \geq0,\;\forall \m{w}\in\C{W}.
\]
Therefore, $\hat{\m{w}}\in \C{W}^*$ is a solution of  $\textrm{VI}(\Theta,\mathcal{J},\C{W})$.

Next, we prove (iii). Followed by (\ref{Sec3-54}) and the positive definiteness of $H$,   the sequence $\{\m{w}^k\}$ is  uniformly bounded. So, there exists a subsequence $\{\m{w}^{k_{j}}\}$ converging to a point $\m{w}^{\infty}=(\m{x}^{\infty},y^\infty,\lambda^\infty)\in\C{W}$.
 Recalling the proof of (ii), we know $\lim\limits_{j\rightarrow\infty}\m{w}^{k_j}  =\m{w}^{\infty}$. So
  $\m{w}^\infty\in\C{W}^*$ is a solution point of $\textrm{VI}(\Theta,\mathcal{J},\C{W})$.
Since (\ref{Sec3-54}) holds for any $\m{w}^*\in\C{W}^*$, by (\ref{Sec3-54}) again and $\m{w}^\infty\in\C{W}^*$, we have for all $l\geq k_{j}$ that
\[
\left\|\m{w}^{l}-\m{w}^{\infty}\right\|_{H}^2+\xi_{3}\left\|\m{E}^l\right\|^2\leq\left
\|\m{w}^{k_{j}}-\m{w}^\infty\right\|_{H}^2+\xi_{3}\left\|\m{E}^{k_{j}}\right\|^2.
\]
{Combining this inequality} with (i), (\ref{Sec3-59bj1})  and the positive definiteness of $H$ {is to obtain} $\lim\limits_{l\rightarrow\infty}\m{w}^{l}=\m{w}^\infty$. Therefore, the whole sequence $\{\m{w}^k\}$ converges to  $\m{w}^\infty$.  $\blacksquare$

Theorem \ref{conver-th} illustrates that our proposed algorithm ADMM-LQP is {globally convergent}.
 In the following, we will show its   sublinear convergence rate for the   ergodic iterates that seems to appear  originally in \cite{BaHZ20}:
\begin{equation}\label{Sec3-62}
\m{w}_{T}:=\frac {1}{1+T}\sum\limits_{k=\kappa}^{\kappa+T}\tilde{\m{w}}^k \quad \textrm{for any } \kappa\geq0.
\end{equation}
\begin{theorem} (Ergodic convergence rate)
Suppose that conditions of Theorem \ref{Theor2} hold.  Then,  for any integers $\kappa\geq0, T>0$ and  for all $k\in[\kappa, \kappa+T]$, there {exists   a constant} $\xi_{3}\geq0$ such that
\begin{equation}\label{Sec3-61}
\Theta(\m{w}_{T})-\Theta(\m{w})+\left\langle\m{w}_{T}-\m{w}, \C{J}(\m{w})\right\rangle
\leq\frac {1}{2(1+T)}\left\{\left\|\m{w}-\m{w}^\kappa\right\|_{H}^2+\xi_{3}\left\|A\m{x}^\kappa+By^\kappa-b\right\|^2\right\}.
\end{equation}
\end{theorem}
\textbf{Proof}.
Summing the inequality (\ref{bj-ke1}) over $k$ between $\kappa$ and $\kappa+T$ gives
\begin{eqnarray}\label{Sec3-63}
\begin{split}
&\sum_{k=\kappa}^{\kappa+T}\Theta(\tilde{\m{w}}^k)-(1+T)\Theta(\m{w}) +\left\langle \sum_{k=\kappa}^{\kappa+T}\tilde{\m{w}}^k-(1+T)\m{w}, \C{J}(\m{w})\right\rangle \\
\leq& \frac{1}{2}\left\{\left\|\m{w}-\m{w}^\kappa\right\|_{H}^2+\xi_{3}\left\|A\m{x}^\kappa+By^\kappa-b\right\|^2\right\},~~\forall \m{w}\in\C{W}.
\end{split}
\end{eqnarray}
By the convexity of $\Theta$ and the definition of $\m{w}_{T}$, we have
\[
\Theta(\m{w}_{T})\leq \frac {1}{1+T}\sum_{k=\kappa}^{\kappa+T}\Theta(\tilde{\m{w}}^k) \textcolor{red}{.}
\]
Dividing (\ref{Sec3-63}) by $(1+T)$  and invoking the above inequality, the conclusion (\ref{Sec3-61}) is confirmed.   $\blacksquare$\medskip

\noindent{\bf Remark  2} {\it
 In the following, we further investigate   another variant of {the sublinear convergence rate  in terms of the objective function value
error and the feasibility error, which looks a bit  more intuitive}. For any $\xi>0$, we define
 $\Gamma_\xi:=\{\lambda|~ \xi\geq \|\lambda\| \}$ and
 \[
 \delta_\xi:=\inf\limits_{\m{u}^*\in\C{X}^*\times \C{Y}^*}\sup\limits_{\lambda\in \Gamma_\xi}
 \left\|\left(\begin{array}{c}
              \m{u}^*-\m{u}^\kappa\\
               \lambda- \lambda^\kappa
             \end{array}\right)\right\|_{H}^2+\xi_{3}\left\|A\m{x}^\kappa+By^\kappa-b\right\|^2,
 \]
where $\m{u}^*=(\m{x}^*;y^*)$.  By setting $\m{w}:=(\m{u}^*,\lambda)$ into (\ref{Sec3-61}) and using the definitions of $\m{w}_{T}$ and $\C{J}(\m{w})$, we get
\[
\begin{split}
& \Theta(\m{w}_{T})-\Theta(\m{w})+\left\langle\m{w}_{T}-\m{w}, \C{J}(\m{w})\right\rangle \\
=&  \Theta(\m{u}_{T})-\Theta(\m{u}^*)- \lambda\tr A\left(\m{x}_T-\m{x}^*\right)- \lambda\tr B\left(y_T-y^*\right)+\left(\lambda_T-\lambda\right)\tr\left(A\m{x}^*+By^*-b\right)\\
=&  \Theta(\m{u}_{T})-\Theta(\m{u}^*)- \lambda\tr \left(A\m{x}_T+By_T -b\right)\\
\leq &\frac {1}{2(1+T)}\left\{\left\|\left(\begin{array}{c}
              \m{u}^*-\m{u}^\kappa\\
               \lambda- \lambda^\kappa
             \end{array}\right)\right\|_{H}^2+\xi_{3}\left\|A\m{x}^\kappa+By^\kappa-b\right\|^2\right\}
\end{split}
\]
where the second equality uses the fact $A\m{x}^*+By^*=b$. Then, it follows from the above inequality that
\begin{eqnarray}\label{Sec3-1jbai}
\begin{split}
&  \Theta(\m{u}_{T})-\Theta(\m{u}^*)+ \xi \left\|A\m{x}_T+By_T -b\right\|\\
=& \sup\limits_{\lambda\in \Gamma_\xi}\left\{ \Theta(\m{u}_{T})-\Theta(\m{u}^*)- \lambda\tr \left(A\m{x}_T+By_T -b\right)\right\}\\
\leq &\frac {1}{2(1+T)}\left\{\inf\limits_{\m{u}^*\in\C{X}^*\times \C{Y}^*}\sup\limits_{\lambda\in \Gamma_\xi}\left\|\left(\begin{array}{c}
              \m{u}^*-\m{u}^\kappa\\
               \lambda- \lambda^\kappa
             \end{array}\right)\right\|_{H}^2+\xi_{3}\left\|A\m{x}^\kappa+By^\kappa-b\right\|^2\right\} \\
=& \frac {\delta_\xi}{2(1+T)}.
\end{split}
\end{eqnarray}

Analogous to the analysis of (\ref{Sec3-1jbai}) together with (\ref{Sec31-1}), we must have \[\Theta(\m{u}_{T})-\Theta(\m{u}^*)- (\lambda^*)\tr \left(A\m{x}_T+By_T -b\right)\geq 0\]
showing that $\Theta(\m{u}_{T})-\Theta(\m{u}^*)\geq -\|\lambda^*\|\left\|A\m{x}_T+By_T -b\right\|.$
Taking $\xi= 2\|\lambda^*\|+1$ in (\ref{Sec3-1jbai}) gives
\[
\begin{split}
& ( \|\lambda^*\|+1) \left\|A\m{x}_T+By_T -b\right\|
\\
\leq &  \Theta(\m{u}_{T})-\Theta(\m{u}^*)+ (2\|\lambda^*\|+1) \left\|A\m{x}_T+By_T -b\right\|
\leq   \frac {\delta_\xi}{2(1+T)},
\end{split}
\]
that is,
\[
\left\|A\m{x}_T+By_T -b\right\|\leq \frac {\delta_\xi}{2(1+T)( \|\lambda^*\|+1)}.
\]
So, we will also get $\Theta(\m{u}_{T})-\Theta(\m{u}^*)\geq-\|\lambda^*\|\left\|A\m{x}_T+By_T -b\right\|
\geq -
\frac {\delta_\xi}{2(1+T)}$ which shows
\[
\mid \Theta(\m{u}_{T})-\Theta(\m{u}^*)\mid \leq \frac {\delta_\xi}{2(1+T)}.
\]
}

Now, we  would investigate the sublinear  convergence rate of our ADMM-LQP for the   primal residuals, {that is, the nonergodic convergence rate.}
 \begin{theorem} (Nonergodic convergence rate)
Suppose that conditions of Theorem \ref{Theor2} hold. Then, for any integer $k>0$, there exists an integer $t\leq k$ such that
\begin{equation} \label{coro2i}
\left\|\m{x}^t-\m{x}^{t-1}\right\|^2\leq \frac{\vartheta}{k}, \quad \textrm{and} \quad
\left\|y^t-y^{t-1}\right\|^2\leq \frac{\vartheta}{k},
\end{equation}
 where   $\vartheta>0$ is a constant  depending
on the problem data and the parameters of ADMM-LQP.
 \end{theorem}
\textbf{Proof}.   Let $k>0$ be any fixed constant and $t\in[1,k]$ be an integer such that
\[
\begin{split}
&  \xi_{1}\sum\limits_{i=1}^p \left\|A_{i}(x_{i}^{t-1}-x_{i}^t)\right\|^2
+\xi_{2}\left\|\m{E}^t\right\|^2
\\=&\min\left\{\xi_{1}\sum\limits_{i=1}^p \left\|A_{i}(x_{i}^{l-1}-x_{i}^l)\right\|^2
+\xi_{2}\left\|\m{E}^l\right\|^2: l=1,2,\cdots,k\right\}.
\end{split}
\]
Then, we deduce by (\ref{Sec3-591ij}) that
\[
\xi_{1}\sum\limits_{i=1}^p \left\|A_{i}(x_{i}^{t-1}-x_{i}^t)\right\|^2
+\xi_{2}\left\|\m{E}^t\right\|^2\leq \frac{\vartheta}{k},
\]
 and   $\vartheta>0$ is a generic constant only depending on the
problem data and the parameters of ADMM-LQP.  Hence, the left   inequality in (\ref{coro2i}) is obtained.
The right inequality in (\ref{coro2i}) holds by combining the  equality (\ref{equ-refor})  and the result $\left\|\m{E}^t\right\|^2\leq  {\vartheta}/{k}$. This completes the proof. $\blacksquare$

\section{Extension of ADMM-LQP}\label{sec445p}
In this part, we consider the following structured convex optimization problem
\begin{equation} \label{Sec4-1j}
\begin{array}{lll}
 \min & \sum\limits_{i=1}^{p}f_i(x_i)+g(y)+h(y)\\
\textrm{s.t. }& \sum\limits_{i=1}^{p}A_ix_i+ By =b,\\
       & x_i \in\mathbb{R}_{+}^{m_i}{,}~y \in\mathbb{R}^d~,i=1,\cdots,p,
\end{array}
\end{equation}
where $f_i(i=1,\cdots,p)$ are the same as before, $h(y)$ is nonsmooth but $g(y)$ is  smooth. We   assume the gradient of $g$ satisfies the Lipschitz condition, i.e., for any $y,\bar{y}\in{\mathbb{R}^d},$  there exists a constant $L_g>0$ such that
 $
 \|g(y)-g(\bar{y})\|\leq L_g  \|y- \bar{y}\|.
 $
 This inequality, by  a Taylor expansion, implies
 \begin{equation}\label{Sec31-01jb}
 g(y)\leq g(\bar{y})+\langle \nabla g(\bar{y}), y- \bar{y}\rangle +\frac{L_g}{2} \|y- \bar{y}\|^2.
\end{equation}
 The function $h(y)$ {is usually} used to promote a data structure different from the structure promoted by $g(y)$ at the solution. Particularly,  constraints of the form $y\in\C{Y}$ where $\C{Y}$ is   a closed convex set, can be incorporated in the optimization problem by letting $h(y)$ be the indicator function of $\C{Y}$, that is, $h(y)=\infty$ when $y\notin\C{Y}$ and otherwise $h(y)=0.$

Using the popular linearization  technique,  the $y$-subproblem in ADMM-LQP can be updated by the following
\begin{equation} \label{Sec4-2j}
 y^{k+1} \leftarrow \arg\min\limits_{y\in\mathbb{R}^d} \left\{
\begin{split}
 &h(y)+ \left\langle y, \nabla g(y^k) -B\tr \lambda^{k+\frac{1}{2}} \right\rangle +\\
 &
 \frac{\beta}{2}   \left\|A\m{x}^{k+1}+By -b\right\|^2+\frac{1}{2}\left\|y-y^k\right\|_{D}^2
 \end{split}\right\},
\end{equation}
where $D$ is a determined proximal matrix.  An obvious advantage of using the  proximal term $\frac{1}{2}\left\|y-y^k\right\|_{D}^2$,  where $D=\sigma\m{I}-\beta B\tr B$ and $\sigma\geq \beta \|B\tr B\|$, is that it could transform the $y$-subproblem into the proximal mapping:
\[
 \textrm{\bf prox}_{h,\sigma}(y_c^k)=\arg\min h(y) +\frac{\sigma}{2} \left\|y-y_c^k\right\|^2.
\]
Here $y_c^k=y^k -\frac{\nabla g(y^k) -B\tr \lambda^{k+\frac{1}{2}}  +\beta B\tr(A\m{x}^{k+1}+By^k -b)}{\sigma}.$

In the following, we briefly analyze the convergence of our extended algorithm named \textbf{eADMM-LQP} for solving problem (\ref{Sec4-1j}). Since the $x_i$-subproblem updates {in} the same way as before, we just need to analyze the first-order optimality condition  of the subproblem in (\ref{Sec4-2j}).

By the convexity of $g(y)$ and (\ref{Sec31-01jb}), we have
\begin{eqnarray}\label{sec34-bw4}
\begin{split}
g(\tilde{y}^k)&\leq  g(y^k)+\left\langle \nabla g(y^k),\tilde{y}^k-y^k\right\rangle +\frac{L_g}{2}\left\|\tilde{y}^k-y^k\right\|^2\nonumber\\
&= g(y^k)+\left\langle \nabla g(y^k),y-y^k+\tilde{y}^k-y\right\rangle +\frac{L_g}{2}\left\|\tilde{y}^k-y^k\right\|^2\nonumber\\
&\leq   g(y)+\left\langle \nabla g(y^k),\tilde{y}^k-y\right\rangle +\frac{L_g}{2}\left\|\tilde{y}^k-y^k\right\|^2.
\end{split}
 \end{eqnarray}
The first-order optimality   {condition   of   $y$-subproblem is}
\[
\begin{split}
&h(y)-h(y^{k+1})+ \left\langle y-y^{k+1}, \nabla g(y^k) -B\tr \lambda^{k+\frac{1}{2}}\right.\\
 &  \left.+\beta B\tr (A\m{x}^{k+1}+By^{k+1} -b)+D(y^{k+1}-y^k)\right\rangle\geq 0, ~~\forall  y \in {\mathbb{R}^d},
\end{split}
\]
which, by    (\ref{Sec2-8})-(\ref{Sec2-9}) and (\ref{Sec2-17}), can be rewritten as
\[
\begin{split}
&h(y)-h(\tilde{y}^k)+ \left\langle y-\tilde{y}^k, \nabla g(y^k) \right\rangle\\
&+\left\langle y-\tilde{y}^k,   -B\tr\tilde{\lambda}^{k}+\alpha B\tr (\lambda^{k}-\tilde{\lambda}^{k})+(\beta B\tr B +D)(\tilde{y}^k-y^k)\right\rangle\geq 0.
\end{split}
\]
{The above inequalities show}
\begin{eqnarray} \label{Sec2-141j}
\begin{split}
&\left\{h(y)+g(y)- h(\tilde{y}^k) -g(\tilde{y}^k)\right\}+\left\langle y -\tilde{y}^k, -B\tr \tilde{\lambda}^{k}\right\rangle \nonumber\\
&\geq
\left\langle  \tilde{y}^k-y, \alpha B\tr (\lambda^{k}-\tilde{\lambda}^{k})+(\beta B\tr B +D)(\tilde{y}^k-y^k)\right\rangle -\frac{L_g}{2}\left\|\tilde{y}^k-y^k\right\|^2,
\end{split}
\end{eqnarray}
which, by combining the previous inequalities (\ref{Sec2-16}) and (\ref{Sec2-19}), gives
\begin{equation}\label{Sec31-4jb}
\begin{split}
\Theta(\m{w})-\Theta(\tilde{\m{w}}^k)  +\left\langle \m{w}-\tilde{\m{w}}^k,  \C{J}(\tilde{\m{w}}^k)+Q(\tilde{\m{w}}^k-\m{w}^k)\right\rangle+\left\|\m{w}^k-\tilde{\m{w}}^k\right\|_{N}^2\geq0
\end{split}
\end{equation}
for any $\m{w}\in\C{W}$, where the notations $\Theta, \m{w}, \tilde{\m{w}}^k, \C{J},Q$ are \textcolor{blue}{the same} as that in Lemma \ref{le3t}, but   the right-lower block of $Q$ is
\begin{equation}\label{Sec31-5jb}
Q_2=\begin{bmatrix}
       \beta B\tr B +D &  -\alpha B\tr  \\
        -B &  \frac 1\beta \m{I}\\
      \end{bmatrix} \textcolor{blue}{,}
\end{equation}
and
\[
N=\begin{bmatrix}
      N_1 & \m{0 }&\m{0 }\\
      \m{0} & \frac{L_{g}}{2} \m{I}& \m{0 }\\
      \m{0 }&\m{0 }&\m{0 }\\
    \end{bmatrix}.
\]
Note that if taking $D=\sigma\m{I}-\beta B\tr B$ and $\sigma\geq \beta \|B\tr B\|$, then the $Q_2$ given by (\ref{Sec31-5jb}) is similar to that in (\ref{Sec2-14}). Next, we use this form   to simplify   convergence analysis of eADMM-LQP.

\begin{theorem}
Let $D=\sigma\m{I}-\beta B\tr B$ and $\sigma\geq \beta \|B\tr B\|$. Then the iterates generated by  eADMM-LQP  satisfy
 (\ref{Sec3-34}) and (\ref{Sec3-38})    with $Q_2$  given by (\ref{Sec31-5jb}).
\end{theorem}
\textbf{Proof}.
Similar to the proof of Theorem \ref{t1-bjc}, we can show
\begin{equation}\label{Sec31-6jb}
\begin{split}
\Theta(\m{w})-\Theta(\tilde{\m{w}}^k)+\left\langle \m{w}-\tilde{\m{w}}^k, \C{J}(\m{w})\right\rangle&\geq\frac {1}{2}\left\{\left\|\m{w}-\m{w}^{k+1}\right\|_{H}^2-\left\|\m{w}-\m{w}^k\right\|_{H}^2\right\}\\
&\quad +\frac {1}{2}\left\|\m{w}^k-\tilde{\m{w}}^k\right\|_{G}^2,
\end{split}
\end{equation}
where $H=\operatorname{diag}(Q_1,H_2)$ with
\[
H_{2}=\begin{bmatrix}
          (1-\frac {\alpha\tau}{\alpha+\tau})\beta B\tr B+D & -\frac \alpha{\alpha+\tau}B\tr  \\
          -\frac \alpha{\alpha+\tau}B  & \frac {1}{(\alpha+\tau)\beta} \m{I} \\
        \end{bmatrix},
\]
and $G=\operatorname{diag}(G_1,G_2)$ with
\[
G_2
=\begin{bmatrix}
          (1-\tau)\beta B\tr B +D-L_{g} \m{I}  & (\tau-1)B\tr  \\
          (\tau-1)B & \frac {2-\alpha-\tau}{\beta}\m{I} \\
        \end{bmatrix}.
\]
Finally, setting $\m{w}=\m{w}^*\in\C{W}^*$ in (\ref{Sec3-34}) ensures the inequality (\ref{Sec3-38}). $\blacksquare$

Analogous to the analytical techniques in Section \ref{se32b}, we  next estimate the lower bound of $\left\|\m{w}^k-\tilde{\m{w}}^k\right\|_{G}^2$. In fact, by a similar way to \textcolor{blue}{the} proof of Lemma \ref{le6-j} together with the previous notations in  (\ref{Sec2-8j0-1}), there exists a constant $\xi_{1}>0$ such that
\begin{equation}\label{Sec4-60}
\begin{split}
&\left\|\m{w}^k-\tilde{\m{w}}^k\right\|_{G}^2\geq \xi_{1}\sum\limits_{i=1}^p \left\|A_{i}\m{E}^{k+1}_{x_{i}}\right\|^2+\beta(1-\alpha)\left\|B\m{E}_y^{k+1}\right\|^2\\&
+\beta(2-\alpha-\tau)\left\|\m{E}^{k+1}\right\|^2+{ \left \| \m{E}_y^{k+1} \right \| ^{2} _{D-L_{g} \m{I}}}
+2\beta(1-\alpha)\left\langle \m{E}^{k+1}, -B\m{E}_y^{k+1}\right\rangle.
\end{split}
\end{equation}
To further investigate the  lower bound of  $\left\|\m{w}^k-\tilde{\m{w}}^k\right\|_{G}^2$, we need {a lemma as follows}.

\begin{lemma}
Suppose $\alpha>-1$. Then the iterates   generated by  eADMM-LQP satisfy
\begin{eqnarray}\label{Sec4-61}
&&\left\langle \m{E}^{k+1}, -B\m{E}_y^{k+1}\right\rangle \geq
\frac {\tau-1}{1+\alpha}\left\langle \m{E}^k, B\m{E}_y^{k+1}\right\rangle
-\frac {1}{\beta(1+\alpha)}\left\langle\m{E}_y^{k+1}, D\m{E}_y^k\right\rangle
\nonumber\\&&+\frac {1}{\beta(1+\alpha)} \left\|\m{E}_y^{k+1}\right\|^{2}_{D}
-\frac {\alpha}{1+\alpha}\left\|B\m{E}_y^{k+1}\right\|^2
-\frac {L_{g}}{2\beta(1+\alpha)} \left({ \left\|\m{E}_y^k\right\|^{2}+ \left\|\m{E}_y^{k+1}\right\|^{2}}\right). \qquad
\end{eqnarray}
\end{lemma}
\textbf{Proof}.
It follows from the optimality condition of $y$-subproblem together with the previous notations in  (\ref{Sec2-8j0-1}) {that
\[
\begin{split}
&h(y)+g(y)-h(y^{k+1})-g(y^{k+1})\\&+\left\langle y-y^{k+1},
-B\tr \lambda^{k+\frac {1}{2}}+\beta B\tr \m{E}^{k+1}
+D\m{E}_y^{k+1}\right\rangle
+\frac{L_{g}}{2}\left\|\m{E}_y^{k+1}\right\|^{2}\geq0.
\end{split}
\]
Letting} $y=y^k$ in the above inequality gives
\begin{equation}\label{Sec4-62}
\begin{split}
&h(y^{k})+g(y^{k})-h(y^{k+1})-g(y^{k+1})\\&-\left\langle \m{E}_y^{k+1},
-B\tr \lambda^{k+\frac {1}{2}}+\beta B\tr \m{E}^{k+1}
+D\m{E}_y^{k+1}\right\rangle
+\frac{L_{g}}{2}\left\|\m{E}_y^{k+1}\right\|^{2}\geq0.
\end{split}
\end{equation}
Similarly, the optimality condition of $y$-subproblem at the $k$-th iteration is
\[
\begin{split}
&h(y)+g(y)-h(y^{k})-g(y^{k})\\&+\left\langle y-y^{k},
-B\tr \lambda^{k-\frac {1}{2}}+\beta B\tr\m{E}^k
+D\m{E}_y^k\right\rangle
+\frac{L_{g}}{2}\left\|\m{E}_y^k\right\|^{2}\geq0,
\end{split}
\]
which, by setting $y=y^{k+1}$, shows
\begin{equation}\label{Sec4-63}
\begin{split}
&h(y^{k+1})+g(y^{k+1})-h(y^{k})-g(y^{k})\\&+\left\langle \m{E}_y^{k+1},
-B\tr \lambda^{k-\frac {1}{2}}+\beta B\tr  \m{E}^k
+D\m{E}_y^k\right\rangle
+\frac{L_{g}}{2}\left\|\m{E}_y^k\right\|^{2}\geq0.
\end{split}
\end{equation}
Combining (\ref{Sec4-62}) and (\ref{Sec4-63}), we have
\begin{equation}\label{Sec4-64}
\begin{split}
&\left\langle \m{E}_y^{k+1},
   B\tr (\lambda^{k+\frac 12} -\lambda^{k-\frac 12}) +\beta B\tr   (\m{E}^k -\m{E}^{k+1}) + D(\m{E}_y^k-\m{E}_y^{k+1} \textcolor{red}{)}
 \right\rangle\\&
+\frac{L_{g}}{2}\left(\left\|\m{E}_y^k\right\|^{2}{ +}\left\|\m{E}_y^{k+1}\right\|^{2}\right)\geq0.
\end{split}
\end{equation}
Then, substituting (\ref{sjb0-0}) into the left\textcolor{red}{-}hand side  of (\ref{Sec4-64}), we deduce
\[
\begin{split}
&\left\langle -B\m{E}_y^{k+1}, (\tau-1)\beta\m{E}^k + (1+\alpha)\beta\m{E}^{k+1} \right\rangle
+\alpha \beta \left\|B\m{E}_y^{k+1}\right\|^2- \left\|\m{E}_y^{k+1}\right\|^2_D\\&
+\left\langle \m{E}_y^{k+1}, D\m{E}_y^k\right\rangle +
\frac{L_{g}}{2}\left(\left\|\m{E}_y^k\right\|^{2}+\left\|\m{E}_y^{k+1}\right\|^{2}\right)
\geq 0,
\end{split}
\]
which is equivalent to (\ref{Sec4-61}).  $\blacksquare$

\begin{theorem} (Lower bound estimation)\label{Theor8}
Let $D=\sigma\m{I}-\beta B\tr B$ with $\sigma\geq \beta \|B\tr B\| + \frac{3-\alpha}{1+\alpha}L_{g}$, and the parameters $(\gamma,\mu)$ satisfy {the conditions in} Theorem \ref{Theor2}. Then, for any
$(\alpha,\tau)\in \mathcal{K},
$
 there exist   $\xi_{1},\xi_{2}>0$ and $\xi_{3}\geq0$  such that
\[
\begin{split}
\left\|\m{w}^k-\tilde{\m{w}}^k \right\|_{G}^2\geq&\xi_{1}\sum\limits_{i=1}^p \left\|A_{i}\m{E}_{x_{i}}^{k+1}\right\|^2
+\xi_{2}\left\|\m{E}^{k+1} \right\|^2
+\xi_{3}\left(\left\|\m{E}^{k+1}\right\|^2-\left\|\m{E}^k\right\|^2\right).
\end{split}
\]
\end{theorem}
\textbf{Proof}.
Plugging (\ref{Sec4-61}) into (\ref{Sec4-60}) and using the Cauchy-Schwartz inequality, we have
\[
\begin{split}
&\left\|\m{w}^k-\tilde{\m{w}}^k \right\|_{G}^2 \\
\geq& \xi_{1}\sum\limits_{i=1}^p \left\|A_{i}\m{E}^{k+1}_{x_{i}}\right\|^2+\beta(1-\alpha)\left\|B\m{E}_y^{k+1}\right\|^2  +\beta(2-\alpha-\tau)\left\|\m{E}^{k+1}\right\|^2+\left \| \m{E}_y^{k+1} \right \| ^{2} _{D-L_{g} \m{I}} \\
& -\frac{\beta}{1+\alpha}\left[ (1-\tau)^2\left\|\m{E}^k\right\|^2  +(1-\alpha)^2\left\|B\m{E}_y^{k+1}\right\|^2\right]- \frac{1-\alpha}{1+\alpha}\left[  \left\|\m{E}_y^k\right\|^2_D  + \left\| \m{E}_y^{k+1}\right\|^2_D\right]\\
&+\frac{2(1-\alpha)}{1+\alpha}\left\| \m{E}_y^{k+1}\right\|^2_D-
 \frac{2\beta\alpha(1-\alpha)}{1+\alpha}\left\|B\m{E}_y^{k+1}\right\|^2 -
\frac{(1-\alpha)L_{g}}{1+\alpha}
\left[\left \| \m{E}_y^{k+1} \right \| ^{2} +\left \| \m{E}_y^k \right \| ^{2}\right]\\
=& \xi_{1}\sum\limits_{i=1}^p \left\|A_{i}\m{E}^{k+1}_{x_{i}}\right\|^2+
\beta\left(2-\alpha-\tau -\frac{(1-\tau)^2}{1+\alpha}\right)\left\|\m{E}^{k+1}\right\|^2\\
& +
\frac{\beta(1-\tau)^2}{1+\alpha}
\left[\left \| \m{E}^{k+1} \right \| ^{2} -\left \| \m{E}^k \right \| ^{2}\right] + R_{xy},
\end{split}
\]
where
\[
\begin{split}
R_{xy}&= \frac{1-\alpha}{1+\alpha}\left[     \left\| \m{E}_y^{k+1}\right\|^2_D- \left\|\m{E}_y^k\right\|^2_D\right]+\left \| \m{E}_y^{k+1} \right \| ^{2} _{D-L_{g} \m{I}}-
\frac{(1-\alpha)L_{g}}{1+\alpha}
\left[\left \| \m{E}_y^{k+1} \right \| ^{2} +\left \| \m{E}_y^k \right \| ^{2}\right]\\
&= \frac{2}{1+\alpha} \left \| \m{E}_y^{k+1} \right \| ^{2} _{D-L_{g} \m{I}}
-\frac{1-\alpha}{1+\alpha} \left \| \m{E}_y^k \right \| ^{2} _{D+L_{g} \m{I}}.
\end{split}
\]
Since  $D=\sigma\m{I}-\beta B\tr B$, then for any $\sigma\geq \beta \|B\tr B\| + \frac{3-\alpha}{1+\alpha}L_{g}$, we have
$D\succeq \frac{3-\alpha}{1+\alpha}L_{g} \m{I}$  and thus  $R_{xy}\geq 0$. So, it holds

\[
\left\|\m{w}^k-\tilde{\m{w}}^k \right\|_{G}^2\geq \xi_{1}\sum\limits_{i=1}^p \left\|A_{i}\m{E}^{k+1}_{x_{i}}\right\|^2
+\xi_{2}\left\|\m{E}^{k+1}\right\|^2 +
\xi_{3}
\left[\left \| \m{E}^{k+1} \right \| ^{2} -\left \| \m{E}^k \right \| ^{2}\right],
\]
where  $\xi_{i} (i=1,2,3)$ are the same as that in Theorem \ref{Theor2}.
  $\blacksquare$

Finally, based on the above discussions and the similar analysis to Section \ref{33glob}, the algorithm
eADMM-LQP converges globally with a sublinear convergence rate in    the ergodic and nonergodic sense.

\section{Concluding remarks}\label{condr}

In this article, we  have  proposed   a partial  LQP-based   ADMM     with larger stepsizes of dual variables as in \cite{BaiLiXuZhang2018} to solve a family of separable convex minimization problems. The so-called LQP regularization term is added to the first grouped subproblems with proper proximal parameters but it is not added to  the final subproblem. With the aid of prediction-correction technique, we establish the global convergence of the proposed algorithm and its sublinear convergence rate in terms of the objective residual and the primal {residual}. The proposed algorithm is also extended to solve a nonsmoooth composite convex optimization, where the second subproblem is updated by using the linearization technique and proximal technique. We analyze  the global sublinear convergence properties of the extended method in a brief way.

Notice that our result about the ergodic convergence rate is established for a new average iterate $\m{w}_{T}:=\frac {1}{1+T}\sum _{k=\kappa}^{\kappa+T}\tilde{\m{w}}^k$ for any integer $\kappa\geq0$, which is more general than the traditional iterate $\m{w}_{T}:=\frac {1}{1+T}\sum _{k=0}^{T}\tilde{\m{w}}^k$ in the literature of ADMM. This makes it possible to use some accelerated techniques in practical programming, see the experiments in \cite{BaHZ20} for more details. Another observation is that the constrained set $\C{Y}$ is any closed convex set, so our problem and algorithm are more general  than two existing researches in \cite{liuGY19,WuLi2019}  and thus have a large number of potential applications.  How to extend our algorithm to the  nonseparable nonconvex optimization problem and how to solve the subproblem inexactly are very interesting {topics. Finally, whether the proposed algorithm could use indefinite proximal term  needs our further investigations in the future work.}

\end{document}